\documentclass[11pt]{amsart}
\usepackage{amssymb}
\usepackage[margin=1.1in]{geometry}
\usepackage{enumerate}

\newtheorem{theorem}{Theorem}[section]
\newtheorem{proposition}[theorem]{Proposition}
\newtheorem{lemma}[theorem]{Lemma}
\newtheorem{corollary}[theorem]{Corollary}
\theoremstyle{definition}

\newtheorem{example}[theorem]{Example}
\newtheorem{remark}[theorem]{Remark}
\newtheorem{problem}[theorem]{Problem}
\theoremstyle{plain}
\newtheorem{introthm}{Theorem}

\newcommand{\R}{\mathbb{R}}
\newcommand{\N}{\mathbb{N}}
\newcommand{\C}{\mathbb{C}}
\newcommand{\com}{\mathrm{com}}
\newcommand{\Tr}{\mathrm{Tr}}
\newcommand{\supp}{\operatorname{supp}}
\newcommand{\Mean}{\mathcal M}

\begin{document}

\title[The commutant of the adjoint center]{The commutant of the adjoint
center need not be the center of the dual}

\author{H. Duru}
\address{Department of Mathematics, Faculty of Science, Istanbul University, Istanbul, Türkiye}
\email{hduru@istanbul.edu.tr}

\author{A. K. Kitover}
\address{Community College of Philadelphia, 1700 Spring Garden St.,
Philadelphia, PA 19130, USA}
\email{akitover@ccp.edu}

\author{M. Orhon}
\address{Department of Mathematics \& Statistics, University of New
Hampshire, Durham, NH 03824, USA}
\email{mo@unh.edu}

\date{November 23, 2026}


\begin{abstract}
Let $X$ be a Banach lattice with ideal center $Z(X)$, and let
$Z^*(X)=\{S^*:S\in Z(X)\}$ be the set of adjoints of central operators,
acting on the dual $X^*$. The center $Z(X^*)$ of the dual is always contained
in the commutant $(Z^*(X))^{\com}$ of $Z^*(X)$, and equality is known when
$X$ has a topological order unit. It has been asserted that equality holds
whenever $X$ has a topologically full center; the published proof of this
assertion is incorrect, and the question has been recorded as open. We show
that the assertion is false. For the Marcinkiewicz sequence space
$m_\psi$, $\psi(t)=\log(1+t)$ (the sequence space underlying the
Dixmier--Macaev ideal), which is Dedekind complete and therefore has a
topologically full center, there is a bounded operator on the dual that
commutes with every adjoint of a central operator but is not itself
central; it maps one Dixmier trace onto another one disjoint from it. The
construction rests on a general sufficient condition, proved with an
invariant mean taking values in a Dedekind complete space $C(K)$ and the
injectivity of $C(K)$. We also collect positive results and a norm-density
criterion for equality. For every Orlicz space $L_\Phi$ we show that every
operator in $(Z^*(L_\Phi))^{\com}$ is lower-triangular with respect to the
decomposition of the dual into its order continuous and singular bands,
with central diagonal entries, and, by a splitting lemma based on Young's
inequality, that the off-diagonal entry vanishes. Hence equality holds for
all Orlicz spaces and their duals, without any $\Delta_2$-condition.
\end{abstract}

\maketitle

\section{Introduction}

Let $X$ be a real or complex Banach lattice. An operator $S\in L(X)$ is
\emph{central} if there is $\kappa\ge0$ with $|Sx|\le\kappa|x|$ for all
$x\in X$ (in the complex case $|\cdot|$ is the complex modulus); the set
$Z(X)$ of central operators is the \emph{ideal center} of $X$. It is a
commutative algebra~\cite[Theorems~3.1.9, 3.1.10 and 3.1.12]{MN}, and the adjoint of a central
operator is central on the dual $X^*$ (Lemma~\ref{l:adjcentral}). Write
\[
Z^*(X)=\{S^*:\ S\in Z(X)\}\subseteq Z(X^*),
\]
and, for a set $\mathcal A\subseteq L(X^*)$, let $\mathcal A^{\com}$ be its
commutant in $L(X^*)$. Since $Z(X^*)$ is commutative and contains
$Z^*(X)$,
\begin{equation}\label{e:inclusion}
Z(X^*)\subseteq(Z^*(X))^{\com}
\end{equation}
for every Banach lattice $X$. The question when equality holds in
\eqref{e:inclusion} arose in the theory of Banach $C(K)$-modules and
operators preserving disjointness developed in~\cite{AAK}. Theorem~9.11
of~\cite{AAK} asserts:
\begin{enumerate}[(1)]
\item if $X$ has a topological order unit (a quasi-interior point), then
$Z(X^*)=(Z^*(X))^{\com}$;
\item the same is true if $X$ has a topologically full center.
\end{enumerate}
Part (1) is proved in~\cite{AAK}. The proof of part (2) is incorrect, as was
pointed out by E.~L.~Arenson (private communication)\footnote{The proof of
part~(2) in~\cite{AAK} was contributed by the second author of the present
paper, who bears sole responsibility for the error. The proof of part~(1)
is due to E.~L.~Arenson.}: it uses commutation of an operator $T\in
(Z^*(X))^{\com}$ with limits of adjoints of central operators, but these
limits exist only in the weak$^*$ operator topology, and $T$ need not be
weak$^*$-continuous. The third author recorded part (2) as an open problem
in~\cite[Remark~4]{Or10}. In the same paper it is shown that, for $X$ with a
topologically full center, $Z(X^*)$ is the image of the bidual $Z(X)''$ under
the Arens extension of the embedding $Z(X)\to Z(X^*)$~\cite[Theorem~1,
Corollary~2]{Or10}; in particular $Z(X^*)$ lies in the weak$^*$ operator
closure of $Z^*(X)$.

In this paper we show that part (2) is false. Let $\psi(t)=\log(1+t)$ and let
$m_\psi$ be the Marcinkiewicz sequence space of all real (or complex)
sequences $x$ with
$\sum_{k=1}^nx^\#_k=O(\log n)$, where $x^\#$ is the decreasing rearrangement
of $|x|$ (Section~\ref{s:space}). This is the sequence space underlying the
Dixmier--Macaev ideal $\mathcal M_{1,\infty}$ of compact operators, on which
Dixmier traces are defined~\cite{Di66,LSZ}. (It is strictly larger than
weak $\ell^1$, $\ell^{1,\infty}=\{x:\sup_nnx^\#_n<\infty\}$, which is not
normable.)

\begin{introthm}\label{ti:counter}
The Banach lattice $X=m_\psi$ is Dedekind complete, so its center is
topologically full, but there is a bounded operator $T$ on $X^*$ that
commutes with every $S^*$, $S\in Z(X)$, and is not central. More precisely,
there are disjoint positive Dixmier traces $s_1,s_2$ on $X$ such that
$Ts_1=s_2$. The same holds for the complex space $m_\psi$
(Remark~\ref{r:complex}).
\end{introthm}

The operator $T$ is necessarily not weak$^*$-continuous
(Remark~\ref{r:final}(a)), so it is invisible to the weak$^*$ operator
topology used in~\cite{Or10}. Theorem~\ref{ti:counter} follows from a general
sufficient condition, Theorem~\ref{t:abstract}. Suppose that the linear span
of the central operators $S$ with $S^2=I$ and $|Sx|=|x|$ for all $x$ is dense
in $Z(X)$,
and that $X^*$ contains a nonzero functional $\varphi_2$ and a functional
$\varphi_1$ disjoint from it such that, for all $S\in Z(X)$, the norm of
$S^*$ on the band generated by $\varphi_2$ is at most a constant multiple of
$\|S^*\varphi_1\|$. Then equality in \eqref{e:inclusion} fails. On
$m_\psi$ the pair $(\varphi_1,\varphi_2)$ is given by two Dixmier traces
taken along disjoint sets of scales (Section~\ref{s:weak}); for both of them
$\|S^*\varphi_i\|$ is comparable to the norm of $S$ in the quotient
$\ell^\infty/c_0$. The proof of Theorem~\ref{t:abstract} combines the
injectivity of $C(K)$ for extremally disconnected $K$ (Goodner, Nachbin,
Kelley) with an invariant mean with values in $C(K)$ (Agnew--Morse,
Kantorovich).

In Section~\ref{s:positive} we collect positive results: equality holds for
closed ideals of Banach lattices with a topological order unit
(Proposition~\ref{p:ideal}), for $X$ with order continuous norm
(Proposition~\ref{p:oc}, a special case of~\cite[Theorem~9.10]{AAK}), and, for every $X$, under a norm-density condition
(Theorem~\ref{t:D}). We also show that the failure of equality is equivalent
to the existence of a nonzero bounded module map between two disjoint bands
of $X^*$ (Proposition~\ref{p:crit}). In Section~\ref{s:orlicz} we treat
Orlicz spaces $L_\Phi$. There the singular band of the dual is an AL-space,
and this forces every operator in $(Z^*(L_\Phi))^{\com}$ to be
lower-triangular with respect to the decomposition of the dual into its
order continuous and singular bands, with central diagonal entries
(Theorem~\ref{t:orlicz}). The off-diagonal entry is a module map from
$L_\Psi$ into the singular band that vanishes on $E_\Psi$. A splitting lemma
based on Young's inequality (Lemma~\ref{l:split}) shows that for $x\in
L_\Phi$ and $h\in L_\Psi$ there is a set $A$ with $1_Ax\in E_\Phi$ and
$1_{\Omega\setminus A}h\in E_\Psi$. From this the off-diagonal entry is zero,
and equality holds for every Orlicz space $L_\Phi(\mu)$, $\mu$
$\sigma$-finite, with no $\Delta_2$-condition on $\Phi$ or its complementary
function (Theorem~\ref{t:orliczall}), and also for its dual
(Corollary~\ref{c:dual}). This includes spaces such as $L_{\exp}[0,1]$, and
indeed every $L_\Phi(\mu)$ with $\mu$ finite and nonatomic and $L_\Phi\ne
E_\Phi$, to which Propositions~\ref{p:ideal} and~\ref{p:oc} and
part~(1) of~\cite[Theorem~9.11]{AAK} do not apply. Section~\ref{s:problems}
contains open problems.

\section{Preliminaries}\label{s:prelim}

From now on all Banach lattices are real unless stated otherwise; the
complex version of the main result is obtained by complexification in
Remark~\ref{r:complex}. For standard facts on Banach lattices
we refer to~\cite{AB,MN,Sch,Za}. For a Banach lattice $E$, $E^*$ is its
dual Banach lattice; it is Dedekind complete. Two elements are
\emph{disjoint}, $u\perp v$, if $|u|\wedge|v|=0$. For $A\subseteq E$,
$A^d=\{v:\ v\perp u\text{ for all }u\in A\}$ is its disjoint complement, and
$\{u\}^{dd}$ is the band generated by $u$. In a Dedekind complete Banach
lattice every band $B$ is a projection band, $E=B\oplus B^d$, and the
corresponding band projection is central. For positive functionals
$\varphi,\chi$ on a Banach lattice $E$ and $0\le x\in E$ we shall use the
Riesz--Kantorovich formulas
\begin{equation}\label{e:RK}
|\varphi|(x)=\sup\{|\varphi(y)|:\ |y|\le x\},\qquad
(\varphi\wedge\chi)(x)=\inf\{\varphi(y)+\chi(x-y):\ 0\le y\le x\}
\end{equation}
(the first for arbitrary $\varphi\in E^*$); see~\cite[Theorem~1.18]{AB} (recall that $E^*=E^\sim$ for a Banach
lattice~\cite[Proposition~1.3.7]{MN}).

\subsection{The center}

The following two lemmas are well known (see~\cite[Theorem~2.60]{AB} and~\cite[Proposition~3.1.2]{MN}); we include their short proofs for
completeness.

\begin{lemma}\label{l:adjcentral}
Let $S\in L(X)$ satisfy $|Sx|\le\kappa|x|$ for all $x\in X$. Then
$|S^*\varphi|\le\kappa|\varphi|$ for all $\varphi\in X^*$. In particular
$Z^*(X)\subseteq Z(X^*)$, and $\|S^*\|=\|S\|$.
\end{lemma}

\begin{proof}
Let $0\le x\in X$. By \eqref{e:RK},
\[
|S^*\varphi|(x)=\sup\{|(S^*\varphi)(y)|:\ |y|\le x\}
=\sup\{|\varphi(Sy)|:\ |y|\le x\}.
\]
If $|y|\le x$, then $|Sy|\le\kappa|y|\le\kappa x$, so
$|\varphi(Sy)|\le|\varphi|(|Sy|)\le\kappa|\varphi|(x)$. Hence
$|S^*\varphi|(x)\le\kappa|\varphi|(x)$ for all $x\ge0$, i.e.\
$|S^*\varphi|\le\kappa|\varphi|$. The equality of norms holds for every
bounded operator.
\end{proof}

\begin{lemma}\label{l:bandpres}
Let $E$ be a Banach lattice, $R\in Z(E)$, and $u,v\in E$ with $u\perp v$.
Then $Ru\perp v$. In particular:
\begin{enumerate}[(a)]
\item if $Ru=v$ and $u\perp v$, then $v=0$;
\item $R$ maps every band $B$ of $E$ into $B$.
\end{enumerate}
\end{lemma}

\begin{proof}
Let $|Rw|\le\kappa|w|$ for all $w$, with $\kappa\ge1$. For $a,b\ge0$ we have
$\kappa a\wedge b\le\kappa a\wedge\kappa b=\kappa(a\wedge b)$. Hence
$|Ru|\wedge|v|\le\kappa|u|\wedge|v|\le\kappa(|u|\wedge|v|)=0$.
(a) If $Ru=v$, then $|v|=|v|\wedge|v|=|Ru|\wedge|v|=0$.
(b) If $u\in B$ and $v\in B^d$, then $Ru\perp v$; so $Ru\in B^{dd}=B$.
\end{proof}

We use the following standard description of the center of a Dedekind
complete Banach lattice.

\begin{proposition}\label{p:centerCK}
Let $E\ne\{0\}$ be a Dedekind complete Banach lattice. Then $Z(E)$, with the
operator norm, composition and the order of $L(E)$ (i.e.\ $S\ge0$ iff
$Sx\ge0$ for all $x\ge0$), is isometrically isomorphic, as an algebra and as
a vector lattice, to $C(K)$ for an extremally disconnected compact Hausdorff
space $K$, with the identity operator corresponding to the constant
function $1$.
\end{proposition}

\begin{proof}[References and outline]
Since $E$ is Dedekind complete, so is the space $L_b(E)$ of order bounded
operators (Riesz--Kantorovich), and $Z(E)=\{S:\ |S|\le\kappa I\text{ for
some }\kappa\}$ is the principal ideal of $L_b(E)$ generated by $I$;
hence $Z(E)$ is a Dedekind complete vector lattice with order unit $I$. On
$Z(E)$ the operator norm coincides with the order unit norm
$\inf\{\kappa:\ |S|\le\kappa I\}$, so $Z(E)$ is an AM-space with unit, hence
lattice isometric to some $C(K)$ by Kakutani's theorem. Composition makes
$Z(E)$ an Archimedean $f$-algebra with unit $I$; in such an algebra the
multiplication is determined by the order and the unit, so the Kakutani
isomorphism is multiplicative. Finally $C(K)$ is Dedekind complete, so $K$
is extremally disconnected~\cite[Proposition~2.1.4]{MN}. For the facts used
see~\cite[Proposition~3.1.7 and Theorem~3.1.11]{MN} ($Z(E)$ is the ideal
generated by $I$, and $\|S\|=\inf\{\kappa:\ |S|\le\kappa I\}$ for $S\in
Z(E)$; cf.~\cite[Theorem~144.2]{Za}), \cite[Theorem~3.1.13]{MN} (an $f$-algebra with unit $e$ is
isomorphic, as an algebra and a vector lattice, to its orthomorphisms, so
its multiplication is determined by the order and $e$), \cite[Chapter~20]{Za}
and~\cite{Wi77}.
\end{proof}

A Banach lattice $X$ has a \emph{topologically full center} if for every
$x\in X$ the closure of $Z(X)x=\{Sx:S\in Z(X)\}$ is the closed ideal
generated by $x$~\cite{Wi81}. Banach lattices with a topological order unit
and $\sigma$-Dedekind complete Banach lattices have topologically full
centers~\cite[Section~2, after Corollary~2]{Or10}. For the space used in
this paper this is proved directly in Lemma~\ref{l:weakbasic}(c).

\subsection{Three classical tools}

\begin{theorem}[Goodner, Nachbin, Kelley~\cite{Go,Na,Ke}]\label{t:GNK}
Let $K$ be an extremally disconnected compact Hausdorff space. If $Y$ is a
linear subspace of a real normed space $W$ and $f:Y\to C(K)$ is a bounded
linear operator, then $f$ has a linear extension $\tilde f:W\to C(K)$ with
$\|\tilde f\|=\|f\|$.
\end{theorem}

\begin{theorem}[{Hahn--Banach--Kantorovich; \cite[Theorem~1.25]{AB} with
$H=\{0\}$}]
\label{t:HBK}
Let $L$ be a real vector space, $\mathcal F$ a Dedekind complete vector
lattice and $p:L\to\mathcal F$ a sublinear map, i.e.\ $p(\Phi+\Phi')\le
p(\Phi)+p(\Phi')$ and $p(c\Phi)=cp(\Phi)$ for real $c\ge0$. Then there is a
linear map $\Mean:L\to\mathcal F$ with $\Mean(\Phi)\le p(\Phi)$ for all
$\Phi\in L$.
\end{theorem}

The third tool is an invariant mean with values in $C(K)$. Its proof follows
Agnew and Morse~\cite{AM}; we give it in detail.

\begin{lemma}\label{l:mean}
Let $K$ be an extremally disconnected compact Hausdorff space, and let $U$ be
an abelian group in which every element $u$ satisfies $u^2=1$. Let $L$ be the
vector space of all bounded maps $\Phi:U\to C(K)$ (bounded means
$\|\Phi\|_L:=\sup_u\|\Phi(u)\|_\infty<\infty$), and for $v\in U$ let
$(\tau_v\Phi)(u)=\Phi(uv)$. There is a linear map $\Mean:L\to C(K)$ such
that for all $\Phi\in L$, $v\in U$ and $h\in C(K)$:
\begin{enumerate}[(i)]
\item $\Mean(\Phi)\ge0$ if $\Phi(u)\ge0$ for all $u$;
\item $\Mean(\Phi)=h$ if $\Phi(u)=h$ for all $u$;
\item $\Mean(\tau_v\Phi)=\Mean(\Phi)$;
\item $\|\Mean(\Phi)\|_\infty\le\|\Phi\|_L$;
\item $\Mean(h\Phi)=h\Mean(\Phi)$, where $(h\Phi)(u)=h\cdot\Phi(u)$.
\end{enumerate}
\end{lemma}

\begin{proof}
Since $K$ is extremally disconnected, $C(K)$ is a Dedekind complete vector
lattice: every nonempty subset of $C(K)$ that is bounded in norm has a
supremum and an infimum in $C(K)$. These are the least upper bound and the
greatest lower bound for the order of $C(K)$; they need not be pointwise
suprema, and we only use their order-theoretic properties.

\emph{Step 1: the sublinear map.} For $\Phi\in L$, $n\in\N$ and
$v_1,\dots,v_n\in U$ (repetitions allowed) put
\[
\Sigma(\Phi;v_1,\dots,v_n)=\sup_{u\in U}\ \frac1n\sum_{k=1}^n\Phi(uv_k)\
\in C(K),
\]
and
\[
p(\Phi)=\inf\{\Sigma(\Phi;v_1,\dots,v_n):\ n\in\N,\ v_1,\dots,v_n\in U\}\in
C(K).
\]
All the elements involved have norm at most $\|\Phi\|_L$, so these suprema
and infima exist, and
\begin{equation}\label{e:pbound}
-\|\Phi\|_L\cdot1\le p(\Phi)\le\|\Phi\|_L\cdot1 .
\end{equation}
Clearly $p(c\Phi)=cp(\Phi)$ for real $c\ge0$ (for $c=0$ both sides are $0$).

\emph{Step 2: subadditivity.} Let $\Phi,\Phi'\in L$, and fix
$v_1,\dots,v_n\in U$ and $w_1,\dots,w_m\in U$. Consider the family
$(v_kw_l)_{k\le n,\,l\le m}$ of $nm$ elements of $U$. For every $u\in U$,
since $v_kw_l=w_lv_k$,
\[
\frac1{nm}\sum_{k,l}(\Phi+\Phi')(uv_kw_l)
=\frac1m\sum_{l=1}^m\Bigl(\frac1n\sum_{k=1}^n\Phi\bigl((uw_l)v_k\bigr)\Bigr)
+\frac1n\sum_{k=1}^n\Bigl(\frac1m\sum_{l=1}^m\Phi'\bigl((uv_k)w_l\bigr)\Bigr).
\]
Each inner average of $\Phi$ is at most $\Sigma(\Phi;v_1,\dots,v_n)$, because
$uw_l\in U$; each inner average of $\Phi'$ is at most
$\Sigma(\Phi';w_1,\dots,w_m)$. Hence the left-hand side is at most
$\Sigma(\Phi;v)+\Sigma(\Phi';w)$ for every $u$, and taking the supremum over
$u$,
\[
p(\Phi+\Phi')\le\Sigma\bigl(\Phi+\Phi';(v_kw_l)_{k,l}\bigr)\le
\Sigma(\Phi;v)+\Sigma(\Phi';w).
\]
Thus $p(\Phi+\Phi')-\Sigma(\Phi';w)$ is a lower bound of the set
$\{\Sigma(\Phi;v)\}$ (over all finite families $v$), so it is at most
$p(\Phi)$; i.e.\ $p(\Phi+\Phi')-p(\Phi)\le\Sigma(\Phi';w)$ for all $w$, and
taking the infimum over $w$ gives $p(\Phi+\Phi')\le p(\Phi)+p(\Phi')$.

\emph{Step 3: invariance.} Let $v\in U$ and take $n=2$, $v_1=1$, $v_2=v$.
For every $u$, using $v^2=1$,
\[
\tfrac12\bigl[(\Phi-\tau_v\Phi)(u)+(\Phi-\tau_v\Phi)(uv)\bigr]
=\tfrac12\bigl[\Phi(u)-\Phi(uv)+\Phi(uv)-\Phi(uv^2)\bigr]=0 .
\]
Hence $p(\Phi-\tau_v\Phi)\le\Sigma(\Phi-\tau_v\Phi;1,v)=0$.

\emph{Step 4: the linear map.} By Theorem~\ref{t:HBK} there is a linear
$\Mean:L\to C(K)$ with $\Mean\le p$.
(ii): if $\Phi\equiv h$, then $\Sigma(\Phi;v)=h$ for every family $v$, so
$p(\Phi)=h$ and likewise $p(-\Phi)=-h$. Hence $\Mean(\Phi)\le h$ and
$-\Mean(\Phi)=\Mean(-\Phi)\le-h$, so $\Mean(\Phi)=h$.
(i): if $\Phi\ge0$, then $\Mean(-\Phi)\le p(-\Phi)\le\Sigma(-\Phi;1)=
\sup_u(-\Phi(u))\le0$, so $\Mean(\Phi)\ge0$.
(iii): $\Mean(\Phi)-\Mean(\tau_v\Phi)=\Mean(\Phi-\tau_v\Phi)\le p(\Phi-
\tau_v\Phi)\le0$. Applying this to $\tau_v\Phi$ in place of $\Phi$, and
using $\tau_v\tau_v\Phi=\Phi$ (as $v^2=1$), we also get
$\Mean(\tau_v\Phi)-\Mean(\Phi)\le0$.
(iv): let $\gamma=\|\Phi\|_L$. Then $-\gamma\cdot1\le\Phi(u)\le\gamma\cdot1$
for all $u$, so by (i) and (ii), $-\gamma\cdot1\le\Mean(\Phi)\le\gamma\cdot
1$, i.e.\ $\|\Mean(\Phi)\|_\infty\le\gamma$.

\emph{Step 5: (v).} First let $h=e$ be the indicator function of a clopen
subset of $K$ (so $e^2=e$, $0\le e\le1$). Write $|\Phi|$ for the map
$u\mapsto|\Phi(u)|$ and $\gamma=\|\Phi\|_L$. For every $\Phi'\in L$ we have
$-|\Phi'|\le\Phi'\le|\Phi'|$, so by (i), $|\Mean(\Phi')|\le\Mean(|\Phi'|)$.
Since $|e\Phi|=e|\Phi|\le\gamma e$ for every $u$, (i) and (ii) give
\[
|\Mean(e\Phi)|\le\Mean(e|\Phi|)\le\Mean(\gamma e)=\gamma e .
\]
Hence $\Mean(e\Phi)$ vanishes at every point where $e=0$, i.e.\
$(1-e)\Mean(e\Phi)=0$, i.e.\ $\Mean(e\Phi)=e\Mean(e\Phi)$. In the same way
(with $1-e$ in place of $e$) $e\Mean((1-e)\Phi)=0$. Therefore
\[
e\Mean(\Phi)=e\Mean(e\Phi)+e\Mean((1-e)\Phi)=e\Mean(e\Phi)=\Mean(e\Phi).
\]
By linearity, (v) holds for every linear combination $h$ of indicators of
clopen sets. An extremally disconnected compact space is totally
disconnected, so the clopen sets separate its points, and the linear
combinations of their indicators form a subalgebra of $C(K)$ that contains
$1$ and separates points; by the Stone--Weierstrass theorem it is dense in
$C(K)$. If $h_j\to h$ uniformly, then by (iv)
$\|\Mean(h_j\Phi)-\Mean(h\Phi)\|_\infty\le\|h_j-h\|_\infty\gamma\to0$ and
$\|h_j\Mean(\Phi)-h\Mean(\Phi)\|_\infty\le\|h_j-h\|_\infty\gamma\to0$, so
(v) holds for $h$.
\end{proof}

\section{A sufficient condition for $Z(X^*)\ne(Z^*(X))^{\com}$}
\label{s:abstract}

For a Banach lattice $X$ let
\[
U(X)=\{S\in Z(X):\ S^2=I\text{ and }|Sx|=|x|\text{ for all }x\in X\}.
\]
Every $u\in U(X)$ satisfies $\|u\|=1$, hence $\|u^*\|=1$; $U(X)$ is an
abelian group under composition (because $Z(X)$ is commutative), with
$u^{-1}=u$. For example, if $P$ is a band projection, then $2P-I\in U(X)$:
$(2P-I)^2=I$, and $|(2P-I)x|=|Px-(I-P)x|=|Px|+|(I-P)x|=|x|$ because
$Px\perp(I-P)x$.

\begin{theorem}\label{t:abstract}
Let $X$ be a real Banach lattice such that the linear span of $U(X)$ is
dense in $Z(X)$ in the operator norm. Suppose that there are
$\varphi_1,\varphi_2\in X^*$ with $\varphi_2\ne0$ and
$\varphi_1\perp\varphi_2$, and a constant $\gamma$, such that, with
$B_2=\{\varphi_2\}^{dd}$,
\begin{equation}\label{e:domination}
\bigl\|S^*|_{B_2}\bigr\|\le\gamma\,\|S^*\varphi_1\|\qquad(S\in Z(X)).
\end{equation}
Then there is a bounded operator $T$ on $X^*$ with $TS^*=S^*T$ for all
$S\in Z(X)$ and $T\varphi_1=\varphi_2$. In particular $T\notin Z(X^*)$, and
$Z(X^*)\ne(Z^*(X))^{\com}$.
\end{theorem}

Here $S^*|_{B_2}$ denotes the restriction of $S^*$ to $B_2$; it maps $B_2$
into itself by Lemmas~\ref{l:adjcentral} and~\ref{l:bandpres}(b).

\begin{proof}
\emph{Step 1: the center of $B_2$.} $B_2$ is a band in the Dedekind complete
Banach lattice $X^*$, hence itself a Dedekind complete Banach lattice, and
$B_2\ne\{0\}$ because $\varphi_2\in B_2$. By Proposition~\ref{p:centerCK}
there are an extremally disconnected compact space $K$ and an isometric
algebra and lattice isomorphism $\rho:C(K)\to Z(B_2)$ with $\rho(1)=I_{B_2}$,
the identity of $B_2$.

\emph{Step 2: the representation of $Z(X)$.} Let $S\in Z(X)$. The restriction
$S^*|_{B_2}$ maps $B_2$ into $B_2$ and satisfies $|S^*\varphi|\le\kappa
|\varphi|$ (Lemma~\ref{l:adjcentral}), so it belongs to $Z(B_2)$. Put
\[
\widehat S=\rho^{-1}\bigl(S^*|_{B_2}\bigr)\in C(K).
\]
For $R,S\in Z(X)$ we have $(RS)^*=S^*R^*$, and $RS=SR$ because $Z(X)$ is
commutative; hence $(RS)^*|_{B_2}=S^*|_{B_2}R^*|_{B_2}$, and since $\rho^{-1}$
is linear and multiplicative, $S\mapsto\widehat S$ is linear and
multiplicative, with $\widehat I=1$. Since $\rho$ is isometric,
$\|\widehat S\|_\infty=\|S^*|_{B_2}\|\le\|S^*\|=\|S\|$. If $u\in U(X)$, then
$\widehat u^{\,2}=\widehat{u^2}=\widehat I=1$ and $\|\widehat u\|_\infty\le1$.

\emph{Step 3: a bounded map on the orbit of $\varphi_1$.} Let
$\mathcal V=\{S^*\varphi_1:\ S\in Z(X)\}$, a linear subspace of $X^*$. Define
$f:\mathcal V\to C(K)$ by
\[
f(S^*\varphi_1)=\widehat S .
\]
This is well defined: if $S^*\varphi_1=R^*\varphi_1$, then
$(S-R)^*\varphi_1=0$, so by \eqref{e:domination}
$\|\widehat S-\widehat R\|_\infty=\|(S-R)^*|_{B_2}\|\le
\gamma\|(S-R)^*\varphi_1\|=0$. It is linear, and by \eqref{e:domination}
$\|f(S^*\varphi_1)\|_\infty=\|S^*|_{B_2}\|\le\gamma\|S^*\varphi_1\|$, so
$\|f\|\le\gamma$.

\emph{Step 4: extension.} By Theorem~\ref{t:GNK} there is a linear operator
$\tilde f:X^*\to C(K)$ extending $f$ with $\|\tilde f\|\le\gamma$.

\emph{Step 5: averaging.} Let $U=U(X)$ and let $\Mean$ be the map of
Lemma~\ref{l:mean} for this group $U$ and the space $K$ of Step 1. For
$\varphi\in X^*$ define $\Phi_\varphi:U\to C(K)$ by
\[
\Phi_\varphi(u)=\widehat u\cdot\tilde f(u^*\varphi).
\]
Since $\|\widehat u\|_\infty\le1$ and $\|u^*\|=1$, we have
$\|\Phi_\varphi(u)\|_\infty\le\gamma\|\varphi\|$ for all $u$, so
$\Phi_\varphi\in L$. Put
\[
q(\varphi)=\Mean(\Phi_\varphi)\in C(K).
\]
The map $\varphi\mapsto\Phi_\varphi$ is linear, so $q:X^*\to C(K)$ is
linear, and $\|q(\varphi)\|_\infty\le\gamma\|\varphi\|$ by
Lemma~\ref{l:mean}(iv).

\emph{Step 6: the module property of $q$.} Let $v\in U$ and $\varphi\in
X^*$. For $u\in U$ we have $u^*v^*=(vu)^*=(uv)^*$. Moreover $\widehat
u=\widehat v\,\widehat{uv}$: indeed $\widehat v\,\widehat{uv}=\widehat{vuv}=
\widehat{uv^2}=\widehat u$. Hence
\[
\Phi_{v^*\varphi}(u)=\widehat u\,\tilde f(u^*v^*\varphi)
=\widehat v\cdot\widehat{uv}\,\tilde f((uv)^*\varphi)
=\widehat v\cdot(\tau_v\Phi_\varphi)(u),
\]
that is, $\Phi_{v^*\varphi}=\widehat v\cdot\tau_v\Phi_\varphi$. By
Lemma~\ref{l:mean}(v) and (iii),
\[
q(v^*\varphi)=\Mean(\widehat v\cdot\tau_v\Phi_\varphi)=\widehat v\,\Mean(\tau_v
\Phi_\varphi)=\widehat v\,\Mean(\Phi_\varphi)=\widehat v\,q(\varphi).
\]
Now fix $\varphi$. Both sides of the identity $q(S^*\varphi)=\widehat
S\,q(\varphi)$ are linear in $S\in Z(X)$, and both are continuous in $S$ for
the operator norm: $\|q(S^*\varphi)\|_\infty\le\gamma\|S\|\|\varphi\|$ and
$\|\widehat S\,q(\varphi)\|_\infty\le\|S\|\,\gamma\|\varphi\|$. They agree on
$U$, hence on its linear span, hence on its closure, which is $Z(X)$ by
hypothesis:
\begin{equation}\label{e:module}
q(S^*\varphi)=\widehat S\,q(\varphi)\qquad(S\in Z(X),\ \varphi\in X^*).
\end{equation}

\emph{Step 7: $q(\varphi_1)=1$.} For $u\in U$, $u^*\varphi_1\in\mathcal V$,
so $\tilde f(u^*\varphi_1)=f(u^*\varphi_1)=\widehat u$, and
$\Phi_{\varphi_1}(u)=\widehat u^{\,2}=1$ for all $u$. By
Lemma~\ref{l:mean}(ii), $q(\varphi_1)=1$.

\emph{Step 8: the operator $T$.} Define
\[
T\varphi=\rho\bigl(q(\varphi)\bigr)\varphi_2\qquad(\varphi\in X^*).
\]
Here $\rho(q(\varphi))\in Z(B_2)$ is applied to $\varphi_2\in B_2$, so
$T\varphi\in B_2\subseteq X^*$. $T$ is linear, and
$\|T\varphi\|\le\|\rho(q(\varphi))\|\,\|\varphi_2\|=\|q(\varphi)\|_\infty
\|\varphi_2\|\le\gamma\|\varphi_2\|\,\|\varphi\|$. Let $S\in Z(X)$. By
\eqref{e:module}, the multiplicativity of $\rho$ and the definition of
$\widehat S$,
\[
T(S^*\varphi)=\rho\bigl(\widehat S\,q(\varphi)\bigr)\varphi_2
=\rho(\widehat S)\bigl(\rho(q(\varphi))\varphi_2\bigr)
=S^*|_{B_2}\bigl(T\varphi\bigr)=S^*T\varphi .
\]
Thus $T\in(Z^*(X))^{\com}$. Finally $T\varphi_1=\rho(1)\varphi_2=\varphi_2$
by Step 7. If $T$ were central, Lemma~\ref{l:bandpres}(a), applied with
$R=T$, $u=\varphi_1$ and $v=\varphi_2$, would give $\varphi_2=0$, a
contradiction.
\end{proof}

\begin{remark}\label{r:U}
The hypotheses used are: the density of the span of $U(X)$, the inequality
\eqref{e:domination}, $\varphi_2\ne0$, and $\varphi_1\perp\varphi_2$.
Positivity of $\varphi_1,\varphi_2$ is not needed.
\end{remark}

\section{The space $m_\psi$ and Dixmier traces}\label{s:weak}

\subsection{The space}\label{s:space}
Let $\psi(t)=\log(1+t)$ for $t\ge0$. For a real sequence $x=(x_n)_{n\ge1}$
tending to $0$, let $x^\#=(x^\#_n)_{n\ge1}$ be the decreasing rearrangement
of $|x|$; equivalently, for $N\in\N$,
\begin{equation}\label{e:maxformula}
\sum_{r=1}^Nx^\#_r=\max\Bigl\{\sum_{n\in J}|x_n|:\ J\subseteq\N,\ \#J=N
\Bigr\}.
\end{equation}
The \emph{Marcinkiewicz sequence space} $X=m_\psi$ consists of all $x\in c_0$
with
\[
\|x\|=\sup_{n\ge1}\frac1{\psi(n)}\sum_{r=1}^nx^\#_r<\infty .
\]
(If $x$ is a sequence with $\limsup_n|x_n|=c>0$, then for every $N$ there is a
set $J$ of $N$ indices with $|x_n|>c/2$ for $n\in J$, so $\sum_{n\in J}|x_n|
\ge Nc/2$, which is not $O(\log N)$; so the restriction to $c_0$ is
harmless.)

\begin{lemma}\label{l:banach}
$(X,\|\cdot\|)$ is a Banach lattice for the coordinatewise order, and it is
an ideal of the space $\R^\N$ of all sequences: if $x\in X$, $y\in\R^\N$
and $|y|\le|x|$, then $y\in X$ and $\|y\|\le\|x\|$. The norm is symmetric:
$\|x\|=\|y\|$ whenever $x^\#=y^\#$.
\end{lemma}

\begin{proof}
If $|y|\le|x|$ coordinatewise, then $y\in c_0$ and, by \eqref{e:maxformula},
$\sum_{r\le N}y^\#_r\le\sum_{r\le N}x^\#_r$ for every $N$; so $y\in X$ and
$\|y\|\le\|x\|$. Symmetry is clear from the definition. For $x,y\in X$ and
$J$ with $\#J=N$, $\sum_{n\in J}|x_n+y_n|\le\sum_J|x_n|+\sum_J|y_n|\le
\sum_{r\le N}x^\#_r+\sum_{r\le N}y^\#_r$; taking the maximum over $J$,
$\|x+y\|\le\|x\|+\|y\|$. Homogeneity and definiteness are clear, so
$\|\cdot\|$ is a lattice norm.

Completeness: let $(x^{(j)})$ be Cauchy in $X$. Since $|x_n|\le x^\#_1\le\psi(1)\|x\|$ for every coordinate $n$, the sequence is
Cauchy coordinatewise; let $x$ be its coordinatewise limit. Fix $\varepsilon>0$
and $j_0$ with $\|x^{(j)}-x^{(i)}\|\le\varepsilon$ for $i,j\ge j_0$. For fixed
$N$, $J$ with $\#J=N$ and $j\ge j_0$, we have $\sum_{n\in J}|x^{(j)}_n-
x^{(i)}_n|\le\varepsilon\psi(N)$ for $i\ge j_0$; letting $i\to\infty$,
$\sum_{n\in J}|x^{(j)}_n-x_n|\le\varepsilon\psi(N)$. Taking the maximum over
$J$, $\sum_{r\le N}(x^{(j)}-x)^\#_r\le\varepsilon\psi(N)$ (we will see in a
moment that $x^{(j)}-x\in c_0$). In particular, for $N=1$,
$\sup_n|x^{(j)}_n-x_n|\le\varepsilon\psi(1)$; since $\varepsilon$ is
arbitrary, $x^{(j)}\to x$ uniformly, so $x\in c_0$, and the
estimate shows $x^{(j)}-x\in X$ and $\|x^{(j)}-x\|\le\varepsilon$ for $j\ge
j_0$. Thus $x\in X$ and $x^{(j)}\to x$.
\end{proof}

Let $c_{00}$ be the space of finitely supported sequences and $X_a$ its
closure in $X$; $X_a$ is a closed ideal. Put
\[
X^*_s=(X_a)^\perp=\{\varphi\in X^*:\ \varphi=0\text{ on }X_a\}.
\]

\begin{lemma}\label{l:band}
$X^*_s$ is a band in $X^*$.
\end{lemma}

\begin{proof}
\emph{Ideal.} Let $\varphi\in X^*_s$ and $|\chi|\le|\varphi|$. For $0\le
x\in X_a$, \eqref{e:RK} gives $|\varphi|(x)=\sup\{|\varphi(y)|:|y|\le x\}=0$,
since every $y$ with $|y|\le x$ lies in the ideal $X_a$. So $|\varphi|=0$ on
$(X_a)_+$, hence $|\chi|=0$ on $(X_a)_+$, and $|\chi(y)|\le|\chi|(|y|)=0$
for $y\in X_a$. Thus $\chi\in X^*_s$. Clearly $X^*_s$ is a linear subspace.
\emph{Order closed.} Let $(\varphi_\alpha)$ be an upward directed family in
$(X^*_s)_+$ with supremum $\varphi$ in $X^*$. For $x\ge0$,
$\varphi(x)=\sup_\alpha\varphi_\alpha(x)$ (see~\cite[Theorem~1.19]{AB}),
so $\varphi=0$ on $(X_a)_+$ and therefore on $X_a$. An order closed ideal is
a band.
\end{proof}

Let $x_0=(1/n)_{n\ge1}$. Since $\sum_{k\le n}1/k\le1+\log n$, we have
$x_0\in X$; put $\nu=\|x_0\|$, so $1\le\nu<\infty$. For an infinite set
$A=\{a_1<a_2<\dots\}\subseteq\N$ let $x_A$ be the sequence with
$(x_A)_{a_j}=1/j$ and $(x_A)_n=0$ for $n\notin A$. Then $x_A^\#=x_0$, so
$x_A\in X$ and $\|x_A\|=\nu$.

\begin{lemma}\label{l:weakbasic}
(a) $X$ is Dedekind complete.
(b) For $g\in\ell^\infty$ and $x\in X$, $gx\in X$ and $\|gx\|\le
\|g\|_\infty\|x\|$. $Z(X)=\{M_g:\ g\in\ell^\infty\}$, where $M_gx=gx$ is the
coordinatewise product, and $\|M_g\|=\|g\|_\infty$.
(c) $X$ has a topologically full center; in fact $Z(X)x$ is the ideal
generated by $x$, for every $x\in X$.
(d) If $g\in c_{00}$, then $M_g^*$ vanishes on $X^*_s$.
(e) $U(X)=\{M_g:\ g_n\in\{-1,1\}\text{ for all }n\}$, and the linear span of
$U(X)$ is dense in $Z(X)$.
\end{lemma}

\begin{proof}
(a) Let $\emptyset\ne D\subseteq X$ be bounded above by $w\in X$. The
coordinatewise supremum $s$ of $D$ exists in $\R^\N$ and satisfies $d_0\le
s\le w$ for any $d_0\in D$, so $|s|\le|d_0|+|w|$ and $s\in X$ by
Lemma~\ref{l:banach}; clearly $s$ is the supremum of $D$ in $X$.

(b) $|gx|\le\|g\|_\infty|x|$, so $gx\in X$ and $\|gx\|\le\|g\|_\infty\|x\|$
by Lemma~\ref{l:banach}; hence $M_g\in Z(X)$ and $\|M_g\|\le\|g\|_\infty$.
Testing on the unit vectors $e_n$ (note $\|e_n\|=1/\psi(1)$) gives
$\|M_g\|\ge|g_n|$, so $\|M_g\|=\|g\|_\infty$. Conversely let $S\in Z(X)$
with $|Sx|\le\kappa|x|$. Then $|Se_n|\le\kappa e_n$, so $Se_n=g_ne_n$ with
$|g_n|\le\kappa$. For $x\in X$ and $n\in\N$ write $x=x_ne_n+y$ with $y_n=0$.
Then $|Sy|\le\kappa|y|$ has $n$-th coordinate $0$, so $(Sx)_n=g_nx_n$. Hence
$S=M_g$.

(c) Let $x\in X$ and let $y$ be in the ideal generated by $x$, i.e.\
$|y|\le c|x|$ for some $c$. Define $g_n=y_n/x_n$ if $x_n\ne0$ and $g_n=0$
otherwise. Then $|g_n|\le c$ and $y=gx=M_gx$. So $Z(X)x$ is the ideal
generated by $x$, and its closure is the closed ideal generated by $x$.

(d) If $g\in c_{00}$, then $gx\in c_{00}\subseteq X_a$, so
$(M_g^*\varphi)(x)=\varphi(gx)=0$ for $\varphi\in X^*_s$.

(e) By (b), $M_g^2=I$ iff $g_n^2=1$ for all $n$, and then $|M_gx|=|x|$. For
$D\subseteq\N$, $1_D=\frac12\bigl(1+(2\cdot1_D-1)\bigr)$, where $1$ is the
constant sequence and $2\cdot1_D-1$ is a $\pm1$ sequence. Hence the span of
$U(X)$ contains every $M_g$ with $g$ taking finitely many values, and these
$g$ are uniformly dense in $\ell^\infty$; by (b) the corresponding $M_g$ are
dense in $Z(X)$.
\end{proof}

\subsection{Means along scales}
Let $\mathcal B$ be the space of bounded real Borel functions on
$[0,\infty)$ with the supremum norm. For $0\le\alpha<\beta$ let
$A_{[\alpha,\beta]}(H)=(\beta-\alpha)^{-1}\int_\alpha^\beta H(u)\,du$ for
$H\in\mathcal B$; then $A_{[\alpha,\beta]}\in\mathcal B^*$,
$\|A_{[\alpha,\beta]}\|=1$.

\begin{lemma}\label{l:omega}
Let $0\le\alpha_k<\beta_k$ with $\alpha_k\to\infty$ and
$\beta_k-\alpha_k\to\infty$. By the Banach--Alaoglu theorem the sequence
$(A_{[\alpha_k,\beta_k]})_k$ has a weak$^*$ cluster point $\omega\in
\mathcal B^*$, i.e.\ for every $\varepsilon>0$, finitely many
$H_1,\dots,H_m\in\mathcal B$ and $k_0$ there is $k\ge k_0$ with
$|A_{[\alpha_k,\beta_k]}(H_i)-\omega(H_i)|<\varepsilon$ for $i\le m$. Every
such $\omega$ satisfies:
\begin{enumerate}[(i)]
\item $\omega(H)\ge0$ if $H\ge0$, and $\omega(1)=1$;
\item $\omega(H(\cdot+c))=\omega(H)$ for $H\in\mathcal B$ and real $c\ge0$;
\item if $H,H'\in\mathcal B$ and $H(u)-H'(u)\to0$ as $u\to\infty$, then
$\omega(H)=\omega(H')$;
\item if $\sup_{u\in[\alpha_k,\beta_k]}|H(u)|\to0$ as $k\to\infty$, then
$\omega(H)=0$.
\end{enumerate}
\end{lemma}

\begin{proof}
We use the following fact: if $G\in\mathcal B$ and
$A_{[\alpha_k,\beta_k]}(G)\to0$ as $k\to\infty$, then $\omega(G)=0$. Indeed,
given $\varepsilon>0$ choose $k_0$ with $|A_{[\alpha_k,\beta_k]}(G)|<
\varepsilon$ for $k\ge k_0$, and then $k\ge k_0$ with
$|A_{[\alpha_k,\beta_k]}(G)-\omega(G)|<\varepsilon$; so $|\omega(G)|<
2\varepsilon$.

(i) holds for every $A_{[\alpha,\beta]}$ and passes to cluster points
(apply the definition with $H_1=H$). (ii): for $c\ge0$,
\[
\bigl|A_{[\alpha_k,\beta_k]}(H(\cdot+c))-A_{[\alpha_k,\beta_k]}(H)\bigr|
=\frac1{\beta_k-\alpha_k}\Bigl|\int_{\beta_k}^{\beta_k+c}H-\int_{\alpha_k}^
{\alpha_k+c}H\Bigr|\le\frac{2c\|H\|_\infty}{\beta_k-\alpha_k}\to0,
\]
so the fact applies to $G=H(\cdot+c)-H$. (iii): $\sup_{u\ge\alpha_k}
|H(u)-H'(u)|\to0$ since $\alpha_k\to\infty$, so
$A_{[\alpha_k,\beta_k]}(H-H')\to0$. (iv): $|A_{[\alpha_k,\beta_k]}(H)|\le
\sup_{[\alpha_k,\beta_k]}|H|\to0$.
\end{proof}

\subsection{Dixmier traces}
For $0\le x\in X$ and real $u\ge0$ put
\[
G_x(u)=\sum_{r\le e^u}x^\#_r,\qquad F_x(u)=\frac{G_x(u)}{\psi(e^u)} .
\]
$G_x$ is a nondecreasing step function and $\psi(e^u)$ is continuous and
positive, so $F_x\in\mathcal B$; indeed $0\le F_x\le\|x\|$ because
$G_x(u)\le\|x\|\psi(\lfloor e^u\rfloor)\le\|x\|\psi(e^u)$.

\begin{lemma}\label{l:sublin}
For $0\le x,y\in X$ and $u\ge0$:
\[
G_{x+y}(u)\le G_x(u)+G_y(u)\le G_{x+y}(u+\log2).
\]
\end{lemma}

\begin{proof}
Let $N=\lfloor e^u\rfloor$. By \eqref{e:maxformula} there is $J$ with $\#J=N$
and $G_{x+y}(u)=\sum_J(x+y)=\sum_Jx+\sum_Jy\le G_x(u)+G_y(u)$. Again by
\eqref{e:maxformula} there are $J_1,J_2$ with $\#J_1=\#J_2=N$,
$G_x(u)=\sum_{J_1}x$ and $G_y(u)=\sum_{J_2}y$. The set $J=J_1\cup J_2$ has at
most $2N\le\lfloor2e^u\rfloor=\lfloor e^{u+\log2}\rfloor$ elements, and since
$x,y\ge0$, $G_x(u)+G_y(u)\le\sum_J(x+y)\le G_{x+y}(u+\log2)$.
\end{proof}

\begin{proposition}\label{p:dixmier}
Let $\omega$ be as in Lemma~\ref{l:omega}. There is a unique linear
functional $\Tr_\omega$ on $X$ with $\Tr_\omega(x)=\omega(F_x)$ for $0\le
x\in X$. It is positive, $\|\Tr_\omega\|\le1$, $\Tr_\omega(x_0)=1$,
$\Tr_\omega\in X^*_s$, and $\Tr_\omega(x)=\Tr_\omega(y)$ whenever $0\le
x,y\in X$ and $x^\#=y^\#$.
\end{proposition}

\begin{proof}
\emph{Additivity on the positive cone.} Let $0\le x,y\in X$. Dividing the
inequalities of Lemma~\ref{l:sublin} by $\psi(e^u)$ gives
\[
F_{x+y}(u)\le F_x(u)+F_y(u)\le H(u),\qquad\text{where }
H(u)=F_{x+y}(u+\log2)\,\frac{\psi(2e^u)}{\psi(e^u)} .
\]
Since $1\le\psi(2t)/\psi(t)\le\log3/\log2$ for $t\ge1$ and
$\psi(2t)/\psi(t)\to1$ as $t\to\infty$, $H\in\mathcal B$ and
$H(u)-F_{x+y}(u+\log2)\to0$ as $u\to\infty$. By Lemma~\ref{l:omega}(i)
(monotonicity and linearity of $\omega$), then (iii), then (ii) with
$c=\log2$,
\[
\omega(F_{x+y})\le\omega(F_x)+\omega(F_y)\le\omega(H)=\omega\bigl(F_{x+y}
(\cdot+\log2)\bigr)=\omega(F_{x+y}).
\]
So $\omega(F_{x+y})=\omega(F_x)+\omega(F_y)$. Clearly $F_{cx}=cF_x$ for
$c\ge0$.

\emph{Extension.} For $x\in X$ put $\Tr_\omega(x)=\omega(F_{x^+})-
\omega(F_{x^-})$. If $x=a-b$ with $a,b\ge0$, then $a+x^-=b+x^+$, so by
additivity $\omega(F_a)+\omega(F_{x^-})=\omega(F_b)+\omega(F_{x^+})$, i.e.\
$\Tr_\omega(x)=\omega(F_a)-\omega(F_b)$; this independence of the
decomposition implies that $\Tr_\omega$ is linear, and it is the only linear
extension. It is positive, and $|\Tr_\omega(x)|\le\Tr_\omega(|x|)=
\omega(F_{|x|})\le\|x\|$ because $F_{|x|}\le\|x\|$.

\emph{Value at $x_0$.} $G_{x_0}(u)=\sum_{k\le e^u}1/k$, and
$\sum_{k\le N}1/k-\log N$ converges as $N\to\infty$, while $\psi(e^u)-u\to0$;
so $F_{x_0}(u)\to1$, and $\Tr_\omega(x_0)=1$ by Lemma~\ref{l:omega}(i),(iii).

\emph{Singularity.} If $x\in c_{00}$, then $G_{|x|}$ is bounded and
$F_{|x|}(u)\to0$, so $\Tr_\omega(|x|)=0$ and $|\Tr_\omega(x)|\le
\Tr_\omega(|x|)=0$. By continuity $\Tr_\omega$ vanishes on $X_a$.

\emph{Symmetry.} $F_x$ depends only on $x^\#$.
\end{proof}

The functionals $\Tr_\omega$ are (commutative) Dixmier traces;
see~\cite{Di66,LSZ}.

\begin{lemma}\label{l:lower}
Let $s\in X^*$ be positive, with $s(x)=s(y)$ whenever $0\le x,y$ and
$x^\#=y^\#$, and with $s(x_0)=1$. For $g\in\ell^\infty$ put
$\lambda(g)=\limsup_{n\to\infty}|g_n|$ (the norm of the class of $g$ in
$\ell^\infty/c_0$). Then
\[
\|M_g^*s\|\ge\frac{\lambda(g)}\nu .
\]
\end{lemma}

\begin{proof}
We may assume $\lambda(g)>0$. Let $0<\varepsilon<\lambda(g)$. Infinitely many
$n$ satisfy $|g_n|>\varepsilon$, so one of the sets $\{n:g_n>\varepsilon\}$,
$\{n:g_n<-\varepsilon\}$ is infinite. Replacing $g$ by $-g$ (which does not
change $\|M_g^*s\|$) we may assume that $A=\{n:g_n>\varepsilon\}$ is
infinite. Then $x_A\ge0$, $\|x_A\|=\nu$, $gx_A\ge\varepsilon x_A$
coordinatewise, and $s(x_A)=s(x_0)=1$ because $x_A^\#=x_0=x_0^\#$. Since $s$
is positive,
\[
(M_g^*s)(x_A)=s(gx_A)\ge\varepsilon s(x_A)=\varepsilon,
\]
so $\|M_g^*s\|\ge\varepsilon/\nu$. Let $\varepsilon\to\lambda(g)$.
\end{proof}

\begin{lemma}\label{l:upper}
If $\varphi\in X^*_s$ and $g\in\ell^\infty$, then $\|M_g^*\varphi\|\le
\lambda(g)\|\varphi\|$. Consequently, for every band $B\subseteq X^*_s$,
$\|M_g^*|_B\|\le\lambda(g)$.
\end{lemma}

\begin{proof}
For $N\in\N$ let $g^{(N)}_n=g_n$ for $n>N$ and $g^{(N)}_n=0$ for $n\le N$.
Then $g-g^{(N)}\in c_{00}$, so $M_g^*\varphi=M_{g^{(N)}}^*\varphi$ by
Lemma~\ref{l:weakbasic}(d), and $\|M_{g^{(N)}}^*\varphi\|\le\|g^{(N)}\|_\infty
\|\varphi\|=\sup_{n>N}|g_n|\,\|\varphi\|$. Let $N\to\infty$. For the second
statement note that $B$ is invariant under $M_g^*$ (Lemma~\ref{l:bandpres}).
\end{proof}

\subsection{Two disjoint Dixmier traces}
Choose $a_1=2$ and define recursively
\[
b_k=a_k^2,\quad a'_k=b_k^2,\quad b'_k={a'_k}^2,\quad a_{k+1}={b'_k}^2
\qquad(k\ge1).
\]
Then $2\le a_k<b_k<a'_k<b'_k<a_{k+1}$, and all these numbers and their
differences tend to $\infty$. Let $\omega_1$ be a weak$^*$ cluster point of
$(A_{[a_k,b_k]})_k$ and $\omega_2$ a weak$^*$ cluster point of
$(A_{[a'_k,b'_k]})_k$ (Lemma~\ref{l:omega}), and put
\[
s_1=\Tr_{\omega_1},\qquad s_2=\Tr_{\omega_2}.
\]
Thus $s_1$ ``sees'' only the scales $u\in[a_k,b_k]$ and $s_2$ only the scales
$u\in[a'_k,b'_k]$ (Lemma~\ref{l:omega}(iv)).

For $0\le x\in X$ and real $t>0$ write $x^\#(t)=x^\#_{\lceil t\rceil}$; then
$\int_0^Nx^\#(t)\,dt=\sum_{r\le N}x^\#_r$ for $N\in\N$, and if $r\in\N$,
$r>t$, then $r\ge\lceil t\rceil$ and $x^\#_r\le x^\#(t)$.

\begin{lemma}\label{l:beta}
Let $0\le x\in X$. For every $b\ge1$ there is $\beta\in[b,2b]$ with
$e^\beta x^\#(e^\beta)\le3\|x\|$.
\end{lemma}

\begin{proof}
The function $v\mapsto e^vx^\#(e^v)$ is a nonnegative Borel function.
Suppose $e^vx^\#(e^v)>3\|x\|$ for all $v\in[b,2b]$. Substituting $t=e^v$,
\[
3\|x\|b<\int_b^{2b}e^vx^\#(e^v)\,dv=\int_{e^b}^{e^{2b}}x^\#(t)\,dt
\le\int_0^{\lceil e^{2b}\rceil}x^\#(t)\,dt=\sum_{r\le\lceil e^{2b}\rceil}
x^\#_r .
\]
The last sum is at most $\|x\|\psi(\lceil e^{2b}\rceil)\le\|x\|\log(2+
e^{2b})\le\|x\|(2b+1)\le3\|x\|b$ for $b\ge1$, a contradiction.
\end{proof}

\begin{proposition}\label{p:disjoint}
$s_1\perp s_2$ in $X^*$.
\end{proposition}

\begin{proof}
Since $s_1,s_2\ge0$, by \eqref{e:RK}, for $0\le x\in X$,
\[
(s_1\wedge s_2)(x)=\inf\{s_1(y)+s_2(x-y):\ 0\le y\le x\}.
\]
So it suffices to find, for each $0\le x\in X$, $x\ne0$, a decomposition
$x=y+z$ with $y,z\ge0$, $s_1(y)=0$ and $s_2(z)=0$. Fix such an $x$.

\emph{Ranks.} Since $x\in c_0$, its positive entries can be listed as
$x_{p_1}\ge x_{p_2}\ge\dots$ with distinct positions $p_r$, so that
$x_{p_r}=x^\#_r$ for all $r$ (if $x$ has finitely many positive entries, the
list is finite). We call $r$ the \emph{rank} of the position $p_r$.

\emph{Thresholds.} By Lemma~\ref{l:beta} there are $\beta_k\in[b_k,2b_k]$
and $\delta_k\in[b'_k,2b'_k]$ with
\[
e^{\beta_k}x^\#(e^{\beta_k})\le3\|x\|,\qquad e^{\delta_k}x^\#(e^{\delta_k})
\le3\|x\|.
\]
The scales are ordered as follows:
\[
\cdots<\delta_{k-1}\le2b'_{k-1}<a_k<b_k\le\beta_k\le2b_k<a'_k<b'_k\le
\delta_k\le2b'_k<a_{k+1}<\cdots
\]
(e.g.\ $2b'_{k-1}<{b'_{k-1}}^2=a_k$ since $b'_{k-1}>2$).

\emph{Definition of $y$ and $z$.} Let $Q=\bigcup_k(e^{\beta_k},e^{\delta_k}]
\subseteq(0,\infty)$. Let $y$ be $x$ restricted to the positions whose rank
lies in $Q$ (i.e.\ $y_{p_r}=x_{p_r}$ if $r\in Q$, and $y_n=0$ at all other
positions $n$), and $z=x-y$, the restriction of $x$ to the positions whose
rank does not lie in $Q$. Then $y,z\ge0$.

\emph{Observation.} If $w$ is the restriction of $x$ to a set of positions,
then by \eqref{e:maxformula} $G_w(u)$ is the sum of at most $\lfloor
e^u\rfloor$ entries of $x$ at positions of that set. An entry of rank
$r>e^\theta$ is at most $x^\#(e^\theta)$. The entries of rank $\le e^\theta$
sum to at most $\sum_{r\le e^\theta}x^\#_r\le\|x\|\psi(e^\theta)\le
\|x\|(\theta+1)$.

\emph{$s_1(y)=0$.} Let $k\ge2$ and $u\in[a_k,b_k]$. Every rank of $y$ lies in
some $(e^{\beta_j},e^{\delta_j}]$; if $j\le k-1$ the rank is $\le
e^{\delta_{k-1}}$, and if $j\ge k$ it is $>e^{\beta_k}$. The entries of the
first kind sum to at most $\|x\|(\delta_{k-1}+1)\le\|x\|(2b'_{k-1}+1)$. At
most $\lfloor e^u\rfloor\le e^{b_k}\le e^{\beta_k}$ entries of the second
kind enter $G_y(u)$, each at most $x^\#(e^{\beta_k})$, so they contribute at
most $e^{\beta_k}x^\#(e^{\beta_k})\le3\|x\|$. Hence, using $\psi(e^u)\ge u
\ge a_k={b'_{k-1}}^2$,
\[
F_y(u)\le\frac{\|x\|(2b'_{k-1}+4)}{u}\le\frac{\|x\|(2b'_{k-1}+4)}
{{b'_{k-1}}^2}\qquad(u\in[a_k,b_k]).
\]
The right-hand side tends to $0$ as $k\to\infty$, so $s_1(y)=\omega_1(F_y)=0$
by Lemma~\ref{l:omega}(iv).

\emph{$s_2(z)=0$.} Let $u\in[a'_k,b'_k]$. No rank of $z$ lies in
$(e^{\beta_k},e^{\delta_k}]$. The entries of $z$ of rank $\le e^{\beta_k}$
sum to at most $\|x\|(\beta_k+1)\le\|x\|(2b_k+1)$. The entries of $z$ of rank
$>e^{\beta_k}$ have rank $>e^{\delta_k}$; at most $\lfloor e^u\rfloor\le
e^{b'_k}\le e^{\delta_k}$ of them enter $G_z(u)$, each at most
$x^\#(e^{\delta_k})$, so they contribute at most $3\|x\|$. Hence
\[
F_z(u)\le\frac{\|x\|(2b_k+4)}{a'_k}=\frac{\|x\|(2b_k+4)}{b_k^2}\qquad
(u\in[a'_k,b'_k]),
\]
which tends to $0$; so $s_2(z)=\omega_2(F_z)=0$ by Lemma~\ref{l:omega}(iv).
\end{proof}

\section{The counterexample}\label{s:counter}

\begin{theorem}\label{t:counter}
Let $X=m_\psi$, $\psi(t)=\log(1+t)$, and let $s_1,s_2$ be the Dixmier traces
of Proposition~\ref{p:disjoint}. Then there is a bounded operator $T$ on
$X^*$ such that $TM_g^*=M_g^*T$ for all $g\in\ell^\infty$ and $Ts_1=s_2$. In
particular $T\in(Z^*(X))^{\com}\setminus Z(X^*)$, although $X$ is Dedekind
complete and has a topologically full center.
\end{theorem}

\begin{proof}
We verify the hypotheses of Theorem~\ref{t:abstract} with
$\varphi_1=s_1$, $\varphi_2=s_2$ and $\gamma=\nu$.

\emph{Density.} By Lemma~\ref{l:weakbasic}(e), the span of $U(X)$ is dense
in $Z(X)$.

\emph{Disjointness and nondegeneracy.} $s_1\perp s_2$ by
Proposition~\ref{p:disjoint}, and $s_2\ne0$ since $s_2(x_0)=1$.

\emph{The domination \eqref{e:domination}.} Every $S\in Z(X)$ is $M_g$ for
some $g\in\ell^\infty$ (Lemma~\ref{l:weakbasic}(b)). The band
$B_2=\{s_2\}^{dd}$ is contained in the band $X^*_s$
(Lemma~\ref{l:band}), because $s_2\in X^*_s$ (Proposition~\ref{p:dixmier}).
By Lemma~\ref{l:upper}, and by Lemma~\ref{l:lower} applied to $s=s_1$
(positive, symmetric and $s_1(x_0)=1$ by Proposition~\ref{p:dixmier}),
\[
\|M_g^*|_{B_2}\|\le\lambda(g)\le\nu\,\|M_g^*s_1\|\qquad(g\in\ell^\infty).
\]
Theorem~\ref{t:abstract} now gives $T$. The last assertion follows from
Lemma~\ref{l:weakbasic}(a),(c).
\end{proof}

\begin{remark}\label{r:final}
(a) \emph{$T$ is not weak$^*$-continuous.} The weak$^*$-continuous operators
on $X^*$ are exactly the adjoints $R^*$ of operators $R\in L(X)$. If $T=R^*$,
then $R$ commutes with every $M_g$. Then $Re_n=RM_{1_{\{n\}}}e_n=
M_{1_{\{n\}}}Re_n$, so $Re_n=r_ne_n$ for some $r_n$, and for $x\in X$,
$(Rx)_ne_n=M_{1_{\{n\}}}Rx=R(x_ne_n)=r_nx_ne_n$; hence $R=M_r$ with
$|r_n|\,\|e_n\|=\|Re_n\|\le\|R\|\,\|e_n\|$, so $|r_n|\le\|R\|$ and
$r\in\ell^\infty$, so $R\in Z(X)$
and $T=R^*$ would be central. (In general this follows
from~\cite[Corollary~3]{Or10}.) This is consistent
with~\cite[Theorem~1]{Or10}: commutation with $Z^*(X)$ does not pass to the
weak$^*$ operator closure of $Z^*(X)$ for operators that are not
weak$^*$-continuous.

(b) \emph{$T$ is regular.} The map $q$ in the proof of
Theorem~\ref{t:abstract} is a bounded linear operator into the Dedekind
complete space $C(K)$; it maps order intervals of $X^*$ into norm bounded,
hence order bounded, subsets of $C(K)$, so it is order bounded and therefore
regular~\cite[Theorem~1.18]{AB}. The map $h\mapsto\rho(h)\varphi_2$ is a
difference of positive maps (split $\varphi_2=\varphi_2^+-\varphi_2^-$; here
$\varphi_2=s_2\ge0$). Hence $T$ is regular, and it is a counterexample also
if the commutant is taken in the space of order bounded operators.

(c) \emph{What the construction uses.} Besides the density of the span of
$U(X)$ in $Z(X)$, only two features of $m_\psi$ are used: there is a functional $s_1$ with $\lambda(g)/\nu\le\|M_g^*s_1\|\le
\lambda(g)$ (Lemmas~\ref{l:lower} and~\ref{l:upper}, and $\|s_1\|\le1$), so
that $g\mapsto M_g^*s_1$ embeds the quotient $\ell^\infty/c_0$
isomorphically, and there is a nonzero singular functional $s_2$ disjoint
from $s_1$.

(d) \emph{$C(Q)$ cannot satisfy the hypotheses of
Theorem~\ref{t:abstract}}, in agreement with part~(1)
of~\cite[Theorem~9.11]{AAK}. Let $X=C(Q)$, $Q$ compact; then $X^*=M(Q)$ is
the space of regular Borel measures, $Z(X)=\{M_g:g\in C(Q)\}$,
$\|M_g^*\mu\|=\int|g|\,d|\mu|$, and the band of a measure $\varphi_2$
consists of the measures absolutely continuous with respect to
$|\varphi_2|$, so that $\|M_g^*|_{\{\varphi_2\}^{dd}}\|$ is the
$|\varphi_2|$-essential supremum of $|g|$, which for continuous $g$ equals
$\sup\{|g(k)|:k\in\supp\varphi_2\}$. Suppose \eqref{e:domination} holds with
$\varphi_1=\mu$. Let $k\in\supp\varphi_2$ and let $V$ be an open
neighbourhood of $k$. By Urysohn's lemma there is $g\in C(Q)$, $0\le g\le1$,
$g(k)=1$, $g=0$ off $V$; then $1\le\gamma\int g\,d|\mu|\le\gamma|\mu|(V)$. By
outer regularity of $|\mu|$, $|\mu|(\{k\})=\inf_V|\mu|(V)\ge1/\gamma$. So
every point of $\supp\varphi_2$ is an atom of $\mu$ of mass at least
$1/\gamma$; there are at most $\gamma\|\mu\|$ such points, so $\supp
\varphi_2$ is finite and $\varphi_2$ is a nonzero finite combination of point
masses at atoms of $\mu$. Then $|\varphi_2|\wedge|\mu|\ne0$, contradicting
$\varphi_1\perp\varphi_2$.
\end{remark}

\begin{remark}\label{r:complex}
\emph{Complex scalars.} Let $X_\C$ be the complexification of the real
Banach lattice $X=m_\psi$. Its center $Z(X_\C)$ is the complexification of
$Z(X)$, i.e.\ $\{M_g:g\in\ell^\infty(\C)\}$, and $(X_\C)^*$ is the
complexification of $X^*$. Let $T_\C$ be the complexification of the
operator $T$ of Theorem~\ref{t:counter}. For $g=g_1+ig_2$ with real $g_1,g_2$
we have $M_g^*=(M_{g_1}^*)_\C+i(M_{g_2}^*)_\C$, which commutes with $T_\C$.
Lemma~\ref{l:bandpres} holds verbatim for complex Banach lattices (with the
complex modulus), so $T_\C s_1=s_2$, $s_1\perp s_2$, $s_2\ne0$ show that
$T_\C$ is not central.
\end{remark}

\begin{remark}\label{r:function}
\emph{The Marcinkiewicz function space.} Theorem~\ref{t:counter} remains
true for the function space $X=M_\psi(0,\infty)$, $\psi(t)=\log(1+t)$, of
all measurable $x$ on $(0,\infty)$ with $\|x\|=\sup_{t>0}\psi(t)^{-1}
\int_0^tx^*<\infty$, where $x^*$ is the decreasing rearrangement of $|x|$.
We only outline the changes; the arguments of Sections~\ref{s:weak}
and~\ref{s:counter} otherwise go through with sums replaced by integrals.

(1) Since $\psi(t)/t\to1$ as $t\to0$, we have $\|x\|_\infty\le\|x\|$, so
$X\subseteq L^\infty$; and $\int_0^tx^*=O(\log t)$ gives $x^*(t)\to0$ as
$t\to\infty$. $X$ is a Dedekind complete, symmetric ideal of $L^0$, and
$Z(X)=\{M_g:\ g\in L^\infty\}$ with $\|M_g\|=\|g\|_\infty$. (The proof of
Lemma~\ref{l:weakbasic}(b) used unit vectors and must be modified: a
central operator $S$ commutes with the band projections $M_{1_B}$, so
$S1_A=g_A1_A$ for sets $A$ of finite measure, with $|g_A|\le\kappa$ and
$g_A=g_{A'}$ on $A\cap A'$; this defines $g\in L^\infty$ with $Sx=gx$ for
simple $x$ with finite-measure support. Every $0\le x\in X$ is the supremum
of an increasing sequence of such functions, since $|\{x>1/n\}|<\infty$, and
central operators are order continuous~\cite[Theorem~2.44]{AB}, so $Sx=gx$.)
As in Lemma~\ref{l:weakbasic}, $Z(X)x$ is the ideal generated by $x$, and
the span of $U(X)=\{M_g:\ |g|=1\}$ is dense in $Z(X)$.

(2) Let $X_a$ be the closure of the functions vanishing outside a set of
finite measure (unlike the sequence case, this is not the order continuous
part of $X$: $\|1_{(0,s)}\|\ge1$ for all $s>0$),
$X^*_s=(X_a)^\perp$ (a band, as in Lemma~\ref{l:band}), and
$\lambda(g)=\inf\{\varepsilon>0:\ |\{|g|>\varepsilon\}|<\infty\}$ for $g\in
L^\infty$. Lemma~\ref{l:upper} holds with the same proof: if
$A=\{|g|>\lambda(g)+\varepsilon\}$, then $|A|<\infty$, $M_{g1_A}$ maps $X$
into $X_a$, and $M_g^*\varphi=M_{g1_{(0,\infty)\setminus A}}^*\varphi$ for
$\varphi\in X^*_s$.

(3) Define $G_x(u)=\int_0^{e^u}x^*$ and $F_x$, $\Tr_\omega$ as in
Section~\ref{s:weak}. Lemma~\ref{l:sublin} becomes
$\int_0^t(x+y)^*\le\int_0^tx^*+\int_0^ty^*\le\int_0^{2t}(x+y)^*$ for
$x,y\ge0$, proved as before using $\int_0^tw^*=\sup\{\int_A|w|:\ |A|=t\}$
in place of \eqref{e:maxformula}, and
Proposition~\ref{p:dixmier} holds with $x_0(t)=1/(1+t)$; here
$\int_0^Tx_0=\psi(T)$, so $\|x_0\|=1$ and $F_{x_0}\equiv1$. If $x$ vanishes
outside a set $A$ of finite measure, then $G_{|x|}\le\|x\|_\infty|A|$ is
bounded, so $\Tr_\omega\in X^*_s$.

(4) In Lemma~\ref{l:lower} replace the infinite set $A\subseteq\N$ by a set
$E\subseteq(0,\infty)$ of infinite measure on which $g>\varepsilon$ (or
$g<-\varepsilon$), and $x_A$ by $x_0\circ\sigma$ on $E$ (and $0$ elsewhere),
where $\sigma:E\to(0,\infty)$ is measure preserving; then the rearrangement
of $x_0\circ\sigma 1_E$ is $x_0$, and we get $\|M_g^*s\|\ge\lambda(g)$.

(5) In Lemma~\ref{l:beta} and Proposition~\ref{p:disjoint}, the role of the
rank is played by a measure preserving map $\rho$ from the support of $x$
onto $(0,\infty)$ with $x=x^*\circ\rho$ a.e.\ on the support of $x$ (Ryff's
theorem in the form~\cite[Chapter~2, Corollary~7.6, p.~83]{BS}; it
applies because $(0,\infty)$ is resonant~\cite[Chapter~2,
Theorem~2.7]{BS} and $x^*(t)\to0$ as $t\to\infty$). We may assume that the support of $x$ has
infinite measure, since otherwise $x\in X_a$ and $s_1(x)=s_2(x)=0$. One puts
$y=x1_{\rho^{-1}(Q)}$ and $z=x-y$, and uses $\int_0^tw^*=\sup\{\int_A|w|:\
|A|=t\}$ in place of \eqref{e:maxformula}; the estimates are the same.

Theorem~\ref{t:abstract} then applies with $\gamma=\|x_0\|=1$.
\end{remark}

\section{Positive results and a criterion}\label{s:positive}

\subsection{Ideals of Banach lattices with a topological order unit}

\begin{proposition}\label{p:ideal}
Let $X$ be a closed ideal of a Banach lattice $Y$ that has a topological
order unit. Then $Z(X^*)=(Z^*(X))^{\com}$.
\end{proposition}

\begin{proof}
\emph{Step 1: $X^*$ as a band of $Y^*$.} Let $X^\perp=\{y^*\in Y^*:\
y^*|_X=0\}$. As in the proof of Lemma~\ref{l:band}, $X^\perp$ is a band in
$Y^*$. Since $Y^*$ is Dedekind complete, $Y^*=X^\perp\oplus B$ with
$B=(X^\perp)^d$; let $P$ be the band projection onto $B$. Let
$\mathrm{r}:B\to X^*$, $\mathrm r(y^*)=y^*|_X$. We claim that $\mathrm r$ is
a surjective lattice isometry.
\begin{itemize}
\item \emph{Lattice homomorphism.} For $y^*\in Y^*$ and $0\le x\in X$, by
\eqref{e:RK}, $|y^*|(x)=\sup\{|y^*(z)|:\ z\in Y,\ |z|\le x\}$, and every such
$z$ lies in the ideal $X$; so $|y^*|(x)=|y^*|_X|(x)$, i.e.\
$|y^*|\,|_X=|y^*|_X|$.
\item \emph{Injective.} If $y^*\in B$ and $y^*|_X=0$, then $y^*\in B\cap
X^\perp=\{0\}$.
\item \emph{Surjective.} Given $x^*\in X^*$, extend it by Hahn--Banach to
$y^*\in Y^*$; then $Py^*\in B$ and $(I-P)y^*\in X^\perp$, so
$\mathrm r(Py^*)=y^*|_X=x^*$.
\item \emph{Isometric.} Let $y^*\in B$. For every $w\in X^\perp$,
$|y^*-w|=|y^*|+|w|\ge|y^*|$ because $y^*\perp w$, so $\|y^*-w\|\ge
\|y^*\|$. Hence $\|y^*|_X\|=\operatorname{dist}(y^*,X^\perp)=\|y^*\|$,
where the first equality is the Hahn--Banach description of the norm of a
restriction.
\end{itemize}

\emph{Step 2: lifting the operator.} Let $T\in(Z^*(X))^{\com}$ and put
\[
\widetilde T=\mathrm r^{-1}T\mathrm rP\in L(Y^*).
\]
Let $R\in Z(Y)$. Since $X$ is an ideal and $|Rx|\le\kappa|x|$, we have
$RX\subseteq X$ and $R|_X\in Z(X)$. By Lemma~\ref{l:adjcentral},
$R^*\in Z(Y^*)$; since $Z(Y^*)$ is commutative and contains $P$, $R^*$
commutes with $P$, and it maps $B$ into $B$ (Lemma~\ref{l:bandpres}(b)). For
$y^*\in B$ and $x\in X$, $(\mathrm rR^*y^*)(x)=y^*(Rx)=(\mathrm ry^*)(R|_Xx)=
((R|_X)^*\mathrm ry^*)(x)$, i.e.\ $\mathrm rR^*=(R|_X)^*\mathrm r$ on $B$,
equivalently $\mathrm r^{-1}(R|_X)^*=R^*\mathrm r^{-1}$ on $X^*$. Hence
\[
\widetilde TR^*=\mathrm r^{-1}T\mathrm rR^*P=\mathrm r^{-1}T(R|_X)^*\mathrm rP
=\mathrm r^{-1}(R|_X)^*T\mathrm rP=R^*\mathrm r^{-1}T\mathrm rP=R^*
\widetilde T .
\]
Thus $\widetilde T\in(Z^*(Y))^{\com}$.

\emph{Step 3.} By~\cite[Theorem~9.11(1)]{AAK}, $\widetilde T\in Z(Y^*)$. A
central operator leaves the band $B$ invariant and its restriction to $B$ is
central; since $\mathrm r$ is a lattice isomorphism,
$T=\mathrm r\,\widetilde T|_B\,\mathrm r^{-1}\in Z(X^*)$. The reverse
inclusion is \eqref{e:inclusion}.
\end{proof}

Proposition~\ref{p:ideal} applies, for example, to $X=C_0(L)$ for every
locally compact $L$ (a closed ideal of $C(\beta L)$), including
non-$\sigma$-compact $L$, for which $X$ has no topological order unit.

\subsection{Order continuous norm}

The following lemma is~\cite[Exercise~3.1.E3, p.~168]{MN}; we include the
short argument for completeness.

\begin{lemma}\label{l:bandcomm}
Let $E$ be a Dedekind complete Banach lattice and $T\in L(E)$. If $T$
commutes with every band projection of $E$, then $T\in Z(E)$.
\end{lemma}

\begin{proof}
For $u\in E$ let $P_u$ be the band projection onto $\{u\}^{dd}$. Then
$Tu=TP_uu=P_uTu\in\{u\}^{dd}$, so $T$ is band preserving. A bounded band
preserving operator on a Banach lattice is central~\cite[Theorem~3.1.12]{MN};
this rests on the theorem of Abramovich, Veksler and Koldunov~\cite{AVK}
that such an operator is order bounded (see also~\cite[Theorem~4.76]{AB}
and~\cite[Section~2]{Or10}).
\end{proof}

\begin{proposition}\label{p:oc}
If $X$ has order continuous norm, then $Z^*(X)=Z(X^*)$ and
$Z(X^*)=(Z^*(X))^{\com}$.
\end{proposition}

Proposition~\ref{p:oc} is a special case of~\cite[Theorem~9.10,
p.~61]{AAK}; we include a direct proof for the reader's convenience.

\begin{proof}
Let $T\in Z(X^*)$. By Lemma~\ref{l:adjcentral}, $T^*\in Z(X^{**})$. Since
the norm of $X$ is order continuous, $X$ (canonically embedded) is an ideal
of $X^{**}$~\cite[Theorem~2.4.2(v)]{MN};
as $|T^*\xi|\le\kappa|\xi|$, $T^*$ maps $X$ into $X$, and $R=T^*|_X\in Z(X)$.
For $\varphi\in X^*$ and $x\in X$, $(R^*\varphi)(x)=\varphi(T^*x)=
(T^*x)(\varphi)=x(T\varphi)=(T\varphi)(x)$, so $T=R^*\in Z^*(X)$. Hence
$Z^*(X)=Z(X^*)$, and every $T\in(Z^*(X))^{\com}$ commutes with all band
projections of $X^*$, so it is central by Lemma~\ref{l:bandcomm}.
\end{proof}

This covers, e.g., $c_0(\Gamma)$ and $\ell^p(\Gamma)$, $1\le p<\infty$, for
uncountable $\Gamma$, and $L^p(\mu)$, $1\le p<\infty$, for non-$\sigma$-finite
$\mu$.

Propositions~\ref{p:ideal} and~\ref{p:oc} do not exhaust the Banach lattices
for which equality holds. To see this we use two simple observations.

\begin{lemma}\label{l:sum}
Let $X=X_1\oplus X_2$ be the direct sum of two Banach lattices (with
coordinatewise order and any lattice norm, e.g.\ $\|(x_1,x_2)\|=
\max(\|x_1\|,\|x_2\|)$). If $Z(X_i^*)=(Z^*(X_i))^{\com}$ for $i=1,2$, then
$Z(X^*)=(Z^*(X))^{\com}$.
\end{lemma}

\begin{proof}
Identify $X^*$ with $X_1^*\oplus X_2^*$, $(\varphi_1,\varphi_2)(x_1,x_2)=
\varphi_1(x_1)+\varphi_2(x_2)$; this is a lattice isomorphism. Let $P_1,P_2$
be the coordinate projections of $X$; they are central, and $P_i^*$ are the
coordinate projections of $X^*$. Let $T\in(Z^*(X))^{\com}$. Since $T$
commutes with $P_1^*$ and $P_2^*$, it has the form $T=T_1\oplus T_2$ with
$T_i\in L(X_i^*)$. If $S\in Z(X_1)$, then $S\oplus0\in Z(X)$ and
$(S\oplus0)^*=S^*\oplus0$; commutation of $T$ with this operator gives
$T_1S^*=S^*T_1$. Thus $T_1\in(Z^*(X_1))^{\com}=Z(X_1^*)$, and likewise
$T_2\in Z(X_2^*)$. If $|T_i\varphi_i|\le\kappa_i|\varphi_i|$, then
$|T(\varphi_1,\varphi_2)|\le\max(\kappa_1,\kappa_2)|(\varphi_1,\varphi_2)|$,
so $T\in Z(X^*)$. The reverse inclusion is \eqref{e:inclusion}.
\end{proof}

\begin{lemma}\label{l:noqi}
Let $\Gamma$ be uncountable. Then $\ell^1(\Gamma)$ is not lattice
isomorphic to a closed ideal of a Banach lattice with a topological order
unit.
\end{lemma}

\begin{proof}
Suppose that $J$ is a closed ideal of a Banach lattice $Y$ with a
topological order unit $u$, and that $V:\ell^1(\Gamma)\to J$ is a lattice
isomorphism; let $c=\|V^{-1}\|$ and $f_\gamma=Ve_\gamma$, where
$(e_\gamma)$ are the unit vectors. Since $u$ is a topological order unit,
$f_\gamma\wedge nu\to f_\gamma$ as $n\to\infty$, so $f_\gamma\wedge
u\ne0$. As $0\le f_\gamma\wedge u\le f_\gamma$ and $J$ is an ideal,
$f_\gamma\wedge u\in J$, and $V^{-1}(f_\gamma\wedge u)$ lies in the ideal
of $\ell^1(\Gamma)$ generated by $e_\gamma$, which is $\R e_\gamma$. Hence
$f_\gamma\wedge u=c_\gamma f_\gamma$ with $0<c_\gamma\le1$. The elements
$c_\gamma f_\gamma$ are pairwise disjoint, so for every finite
$F\subseteq\Gamma$
\[
\sum_{\gamma\in F}c_\gamma f_\gamma=\sup_{\gamma\in F}(f_\gamma\wedge u)\le u,
\qquad\text{and}\qquad
\sum_{\gamma\in F}c_\gamma=\Bigl\|\sum_{\gamma\in F}c_\gamma e_\gamma
\Bigr\|_1\le c\,\|u\| .
\]
Thus $\sum_{\gamma\in\Gamma}c_\gamma<\infty$, which is impossible for
uncountably many $c_\gamma>0$.
\end{proof}

\begin{example}\label{ex:sum}
Let $\Gamma$ be uncountable and $X=\ell^1(\Gamma)\oplus\ell^\infty$. By
Proposition~\ref{p:oc} (for $\ell^1(\Gamma)$), Proposition~\ref{p:ideal}
(for $\ell^\infty$, which has an order unit) and Lemma~\ref{l:sum},
$Z(X^*)=(Z^*(X))^{\com}$. But Proposition~\ref{p:oc} does not apply to
$X$, whose norm is not order continuous, and Proposition~\ref{p:ideal}
does not apply either: if $X$ were a closed ideal of a Banach lattice $Y$
with a topological order unit, then so would be its closed ideal
$\ell^1(\Gamma)\oplus0$ (an ideal of an ideal is an ideal), contradicting
Lemma~\ref{l:noqi}. (For countable $\Gamma$ the space
$\ell^1\oplus\ell^\infty$ has a topological order unit, so
Proposition~\ref{p:ideal} applies.)
\end{example}

\subsection{A criterion}

\begin{proposition}\label{p:crit}
For a Banach lattice $X$ the following are equivalent:
\begin{enumerate}[(a)]
\item $Z(X^*)\ne(Z^*(X))^{\com}$;
\item there are disjoint bands $B_1,B_2$ of $X^*$ and a nonzero bounded
linear $\Theta:B_1\to B_2$ such that $\Theta(R\varphi)=R\Theta(\varphi)$ for
all $R\in Z^*(X)$ and $\varphi\in B_1$.
\end{enumerate}
\end{proposition}

\begin{proof}
(b)$\Rightarrow$(a). Let $P_1$ be the band projection onto $B_1$ and
$T=\Theta P_1$, regarded as an operator on $X^*$ with values in $B_2$. Each
$R\in Z^*(X)$ is central, hence commutes with $P_1$ (both lie in the
commutative algebra $Z(X^*)$) and leaves $B_1$ and $B_2$ invariant; so
$TR=\Theta P_1R=\Theta RP_1=R\Theta P_1=RT$. Choose $\varphi\in B_1$ with
$\Theta\varphi\ne0$; then $T\varphi=\Theta\varphi\in B_2$ is disjoint from
$\varphi$, so $T$ is not central by Lemma~\ref{l:bandpres}(a).

(a)$\Rightarrow$(b). Let $T\in(Z^*(X))^{\com}\setminus Z(X^*)$. As in the
proof of Lemma~\ref{l:bandcomm}, $T$ is not band preserving: there is
$\varphi$ with $T\varphi\notin\{\varphi\}^{dd}$. Let $P$ be the band
projection onto $\{\varphi\}^{dd}$, $B_1=PX^*$, $B_2=(I-P)X^*$ and
$\Theta=(I-P)T|_{B_1}$. Then $\Theta$ commutes with $Z^*(X)$ (as $P$ does),
and $\Theta\varphi=(I-P)T\varphi\ne0$.
\end{proof}

\subsection{A norm-density condition}

\begin{theorem}\label{t:D}
Let $X$ be a Banach lattice such that for every $0\le\theta\in X^*$ and
every band projection $P$ of $X^*$,
\begin{equation}\label{e:D}
P\theta\in\overline{\{S^*\theta:\ S\in Z(X)\}}\qquad\text{(norm closure)}.
\end{equation}
Then $Z(X^*)=(Z^*(X))^{\com}$.
\end{theorem}

\begin{proof}
Let $T\in(Z^*(X))^{\com}$, $\varphi\in X^*$, and let $P$ be a band
projection of $X^*$. Put $\theta=|\varphi|+|T\varphi|$ and choose
$S_n\in Z(X)$ with $S_n^*\theta\to P\theta$. The operators $D_n=S_n^*-P$
are central on $X^*$. For a central operator $D$ on the Dedekind complete
$X^*$ and $\chi\ge0$ we have $|D\chi|=|D|\chi$ and $|D\eta|\le|D||\eta|$ for
all $\eta$ (orthomorphisms are order bounded and disjointness preserving, so
$|D\eta|=|D||\eta|$ by~\cite[Theorem~2.40]{AB}).
Hence
\[
|D_n\varphi|\le|D_n||\varphi|\le|D_n|\theta=|D_n\theta|,\qquad
|D_nT\varphi|\le|D_n||T\varphi|\le|D_n\theta|,
\]
and $\|D_n\theta\|\to0$ gives $S_n^*\varphi\to P\varphi$ and
$S_n^*T\varphi\to PT\varphi$. Therefore
\[
TP\varphi=\lim_nTS_n^*\varphi=\lim_nS_n^*T\varphi=PT\varphi .
\]
So $T$ commutes with every band projection of $X^*$, and $T\in Z(X^*)$ by
Lemma~\ref{l:bandcomm}.
\end{proof}

\begin{remark}
(a) For $X=C(Q)$, \eqref{e:D} holds: $X^*=M(Q)$, a band projection maps
$0\le\theta$ to $1_E\theta$ for a Borel set $E$, and by Lusin's theorem there
are $f_n\in C(Q)$ with $f_n\to1_E$ in $L^1(\theta)$, i.e.\
$\|f_n\theta-1_E\theta\|\to0$.

(b) Since $\{S^*\theta\}$ is a linear subspace, by Mazur's theorem its norm
closure equals its closure in the weak topology $\sigma(X^*,X^{**})$. By
\cite[Corollary~2]{Or10}, for $X$ with a topologically full center $P\theta$
always lies in the $\sigma(X^*,X)$-closure. So \eqref{e:D} asks exactly for
the passage from the weak$^*$ to the weak closure.

(c) For $X=m_\psi$, \eqref{e:D} fails for $\theta=s_1+s_2$ and the projection
$P$ onto $\{s_2\}^{dd}$, so that $P\theta=s_2$. For every $g\in\ell^\infty$,
$M_g^*\theta-s_2=M_g^*s_1+M_{g-1}^*s_2$ with $M_g^*s_1\perp M_{g-1}^*s_2$, so
by Lemma~\ref{l:lower}
\[
\|M_g^*\theta-s_2\|\ge\max\bigl(\|M_g^*s_1\|,\|M_{g-1}^*s_2\|\bigr)\ge
\tfrac1\nu\max\bigl(\lambda(g),\lambda(1-g)\bigr)\ge\tfrac1{2\nu},
\]
because $\max(|g_n|,|1-g_n|)\ge\frac12$ for every $n$.
\end{remark}

\section{Orlicz spaces}\label{s:orlicz}

In this section all Banach lattices are real; $\mu$ is a $\sigma$-finite
measure on $(\Omega,\Sigma)$, $\Phi$ is an $N$-function and $\Psi$ is its
complementary $N$-function. Let $X=L_\Phi=L_\Phi(\mu)$ with the Luxemburg
norm, and let $E_\Phi$ be the closure in $X$ of the simple functions
supported on sets of finite measure;
$E_\Phi=\{x:\ \int\Phi(c|x|)\,d\mu<\infty\text{ for all }c>0\}$, and $E_\Phi$
is the order continuous part $X_a$ of $X$. We use the following standard
facts.
\begin{enumerate}[(O1)]
\item $X^*=L_\Psi\oplus S$, where $L_\Psi$ (with the Orlicz norm, acting by
$h(x)=\int hx\,d\mu$) is the band of order continuous functionals and
$S=(E_\Phi)^\perp$ is the band of singular functionals. The norm of
$L_\Psi$ has the Fatou property.
\item $Z(X)=\{M_g:\ g\in L^\infty\}$, $\|M_g\|=\|g\|_\infty$, and every $y$ with
$0\le y\le x$ is of the form $y=hx$ with $h\in L^\infty$, $0\le h\le1$
(cf.\ Remark~\ref{r:function}(1)).
\end{enumerate}
For (O1) see~\cite[Proposition~2 and Theorem~4]{Ra}: every $x^*\in X^*$
decomposes uniquely as $x^*=y^*+z^*$ with $z^*$ vanishing on $E_\Phi$ and
$y^*(x)=\int hx\,d\mu$ for some $h\in L_\Psi$ (\cite[Theorem~1 and the
Corollary to Lemma~8]{Ra}, since $\Psi$ is continuous and $\mu$ is
$\sigma$-finite; for finite $\mu$ see~\cite{An}), and $|y^*|\wedge|z^*|=0$
(the remark of Zaanen following~\cite[Proposition~2]{Ra}). In the language
of Riesz spaces: the band of order continuous functionals on $L_\Phi$ is
the associate space of $L_\Phi$, which is $L_\Psi$ with the Orlicz norm, and
$X^*$ is the direct sum of this band and its disjoint complement, the band
of singular functionals~\cite[Theorem~132.5 and pp.~597--598]{Za}.
In~\cite{Ra}
$L_\Phi$ carries the Orlicz norm and $L_\Psi$ the Luxemburg norm; with the
Luxemburg norm on $L_\Phi$, as here, the roles of the two norms are
interchanged: the Orlicz norm of $h\in L_\Psi$ is
$\sup\{\int hx\,d\mu:\ \|x\|\le1\}$, the supremum over the unit ball of the
Luxemburg norm~\cite[(9.25) and Section~14.5]{KR}. The Fatou property of
the Orlicz norm follows from its definition as
$\sup\{\int|h|v\,d\mu:\ \int\Phi(v)\,d\mu\le1\}$~\cite[(9.4)]{KR} and
Fatou's lemma. (In~\cite{KR} the measure is Lebesgue measure on a bounded
set; these two arguments do not depend on it.) (O2) is elementary: for $0\le y\le x$ take $h=y/x$ on
$\{x\ne0\}$ and $h=0$ elsewhere; the description of $Z(X)$ is proved as in
Remark~\ref{r:function}(1).

For $x\in L_\Phi$ put
\[
\theta_\Phi(x)=\inf\Bigl\{\lambda>0:\ \int\Phi(|x|/\lambda)\,d\mu<\infty
\Bigr\},
\]
and define $\theta_\Psi$ on $L_\Psi$ in the same way. Then $\theta_\Phi$ is
a seminorm (by convexity of $\Phi$), $\theta_\Phi\le\|\cdot\|$, $\theta_\Phi$
vanishes exactly on $E_\Phi$, $\theta_\Phi(x)\le\theta_\Phi(y)$ if $|x|\le|y|$,
and
\begin{equation}\label{e:thetamax}
\theta_\Phi(1_{A\cup A'}x)=\max\bigl(\theta_\Phi(1_Ax),\theta_\Phi(1_{A'}x)
\bigr)\qquad(A,A'\in\Sigma),
\end{equation}
because the set of admissible $\lambda$ is an up-ray and
$\int_{A\cup A'}\Phi(|x|/\lambda)<\infty$ if and only if both integrals over
$A$ and $A'$ are finite. In the same way $\theta_\Phi(x\vee y)=
\max(\theta_\Phi(x),\theta_\Phi(y))$ for $x,y\ge0$. The same holds for
$\theta_\Psi$, and $\theta_\Psi(k)\le\|k\|$ for the norm of $L_\Psi$ as a
subspace of $X^*$ (the Orlicz norm, which dominates the Luxemburg norm).

\begin{lemma}\label{l:SAL}
$\operatorname{dist}(x,E_\Phi)=\theta_\Phi(x)$ for $x\in X$, and the norm of
$X^*$ is additive on $S_+$, i.e.\ $S$ is an AL-space.
\end{lemma}

\begin{proof}
If $e\in E_\Phi$, then $\theta_\Phi(x)=\theta_\Phi(x-e)\le\|x-e\|$, so
$\theta_\Phi(x)\le\operatorname{dist}(x,E_\Phi)$. Conversely let
$\lambda>\theta_\Phi(x)$, let $F_1\subseteq F_2\subseteq\cdots$ be sets of
finite measure with union $\Omega$, and $x_n=x1_{\{|x|\le n\}\cap F_n}\in
E_\Phi$. Then $\int\Phi(|x-x_n|/\lambda)\to0$ by dominated convergence, so
$\|x-x_n\|\le\lambda$ for large $n$. Hence the quotient norm of the Banach
lattice $X/E_\Phi$ is $[x]\mapsto\theta_\Phi(x)$. Positive elements of
$X/E_\Phi$ have positive representatives, and $[x]\vee[y]=[x\vee y]$, so
$\|[x]\vee[y]\|=\max(\|[x]\|,\|[y]\|)$ for $[x],[y]\ge0$: $X/E_\Phi$ is an
AM-space. Its dual is an AL-space~\cite[Proposition~1.4.7(i)]{MN}. Finally,
the quotient map $\mathrm q:X\to X/E_\Phi$ is a lattice
homomorphism~\cite[Proposition~1.3.13]{MN} and is interval preserving (if
$0\le\mathrm q(y)\le\mathrm q(x)$, then $\mathrm q(y)=\mathrm q(y^+\wedge x)$
with $0\le y^+\wedge x\le x$), so its adjoint is a lattice
homomorphism~\cite[Theorem~1.4.19(ii)]{MN}; being an isometry of
$(X/E_\Phi)^*$ onto $(E_\Phi)^\perp$, it is a lattice isometry onto
$S=(E_\Phi)^\perp$. (For the Orlicz norm on $L_\Phi$ the additivity of the
norm on $S_+$ is~\cite[Lemma~6]{Ra}; see also~\cite[Theorem~133.6]{Za},
where the result is attributed to And\^o~\cite{An}.)
\end{proof}
For $g\in L^\infty$ and $\varphi\in X^*$ we write $g\varphi=M_g^*\varphi$; on
$L_\Psi$ this is the pointwise product. Both bands are invariant under
$\varphi\mapsto g\varphi$. Let $P$ be the band projection of $X^*$ onto
$L_\Psi$ and $Q=I-P$ the band projection onto $S$. Since $L_\Psi$ is an
ideal of $L^0(\mu)$, its band projections are the operators
$k\mapsto1_Ek$, $E\in\Sigma$. Let $(L_\Psi)_a$ be the order continuous part
of $L_\Psi$; it is $E_\Psi$. (For instance, $L_\Psi=E_\Psi$ if $\Psi$
satisfies $\Delta_2$ for all $u$ when $\mu$ is nonatomic and infinite, for
large $u$ when $\mu$ is nonatomic and finite, and for small $u$ when $\mu$ is
the counting measure on $\N$.)

We call a bounded linear map $\Theta$ from a band $B_1$ to a band $B_2$ of
$X^*$ a \emph{module map} if $\Theta(g\varphi)=g\Theta(\varphi)$ for all
$g\in L^\infty$, $\varphi\in B_1$.

\subsection{Singular functionals}

\begin{lemma}\label{l:Ssupport}
Let $0\le\sigma\in S$. There are $B_m\in\Sigma$ ($m\in\N$) with
$\|1_{\Omega\setminus B_m}\sigma\|\le2^{-m}$ and $1_{B_m}\to0$ $\mu$-a.e.
\end{lemma}

\begin{proof}
Let $F_1\subseteq F_2\subseteq\cdots$ be sets of finite measure with union
$\Omega$. Choose $0\le y_m\in X$, $\|y_m\|\le1$, with
$\sigma(y_m)>\|\sigma\|-2^{-m}$. Since $\int\Phi(y_m)\,d\mu\le1$, we have
$\mu\{y_m>n\}\le1/\Phi(n)\to0$, so there is $n_m$ with $\mu\{y_m>n_m\}<2^{-m}$.
Put $B_m=\{y_m>n_m\}\cup(\Omega\setminus F_m)$. Then $y_m1_{\Omega\setminus
B_m}$ is bounded and supported in $F_m$, so it lies in $E_\Phi$ and $\sigma$
vanishes on it. Hence
\[
\|1_{B_m}\sigma\|\ge(1_{B_m}\sigma)(y_m)=\sigma(y_m)>\|\sigma\|-2^{-m}.
\]
Since $1_{B_m}\sigma$ and $1_{\Omega\setminus B_m}\sigma$ are positive
elements of $S$ with sum $\sigma$, Lemma~\ref{l:SAL} gives
$\|1_{\Omega\setminus B_m}\sigma\|=\|\sigma\|-\|1_{B_m}\sigma\|<2^{-m}$.
Finally, $\sum_m\mu\{y_m>n_m\}<\infty$, so by the Borel--Cantelli lemma
almost every point lies in only finitely many of the sets $\{y_m>n_m\}$; it
also lies in $F_m$ for all large $m$. Hence $1_{B_m}\to0$ a.e.
\end{proof}

\begin{lemma}\label{l:Sdisjoint}
Let $\sigma_1,\sigma_2\in S_+$ with $\sigma_1\perp\sigma_2$, and let
$\varepsilon>0$. There is $h\in L^\infty$, $0\le h\le1$, with
$\|h\sigma_1\|\le2\varepsilon$ and $\|(1-h)\sigma_2\|\le2\varepsilon$.
\end{lemma}

\begin{proof}
Choose $0\le x_i\in X$, $\|x_i\|\le1$, with $\sigma_i(x_i)>\|\sigma_i\|-
\varepsilon$, and put $x=x_1+x_2$. By \eqref{e:RK} there is $0\le y\le x$
with $\sigma_1(y)+\sigma_2(x-y)<\varepsilon$, and by (O2) $y=hx$ with
$h\in L^\infty$, $0\le h\le1$. Since $\sigma_1\ge0$ and $hx_1\le hx$,
$\sigma_1(hx_1)<\varepsilon$, so
\[
\|(1-h)\sigma_1\|\ge\sigma_1((1-h)x_1)=\sigma_1(x_1)-\sigma_1(hx_1)>
\|\sigma_1\|-2\varepsilon .
\]
By Lemma~\ref{l:SAL}, $\|h\sigma_1\|=\|\sigma_1\|-\|(1-h)\sigma_1\|<
2\varepsilon$. In the
same way, using $\sigma_2((1-h)x_2)\le\sigma_2((1-h)x)<\varepsilon$, we get
$\|(1-h)\sigma_2\|<2\varepsilon$.
\end{proof}

\begin{lemma}\label{l:Smodule}
Let $B_1,B_2$ be disjoint bands of $X^*$ contained in $S$. Every module map
$\Theta:B_1\to B_2$ is zero.
\end{lemma}

\begin{proof}
It suffices to show $\Theta\sigma=0$ for $0\le\sigma\in B_1$. Put
$\tau=\Theta\sigma$; then $|\tau|\in B_2$, so $|\tau|\perp\sigma$. Let
$\varepsilon>0$ and let $h$ be given by Lemma~\ref{l:Sdisjoint} for
$\sigma_1=\sigma$, $\sigma_2=|\tau|$. Since $M_h^*$ is a positive central
operator, $|h\tau|=h|\tau|$, and
\[
\|\tau\|\le\|h|\tau|\|+\|(1-h)|\tau|\|\le\|h\tau\|+2\varepsilon
=\|\Theta(h\sigma)\|+2\varepsilon\le2(\|\Theta\|+1)\varepsilon .
\]
So $\tau=0$.
\end{proof}

\begin{remark}
Lemma~\ref{l:Smodule} is where Orlicz spaces differ from $m_\psi$. In the
counterexample of Section~\ref{s:counter}, $s_1$ and $s_2$ are disjoint
singular functionals, but the band $X^*_s$ of $m_\psi$ is not an AL-space
(otherwise the proofs of Lemmas~\ref{l:Sdisjoint} and~\ref{l:Smodule}, with
Lemma~\ref{l:weakbasic}(c) in place of (O2), would apply to the restriction
of the operator $T$ of Theorem~\ref{t:counter} to $\{s_1\}^{dd}$),
and $\|M_g^*s_i\|$ is comparable to $\lambda(g)$ for both $i$, so no
analogue of Lemma~\ref{l:Sdisjoint} holds.
\end{remark}

\subsection{The structure of the commutant}

\begin{theorem}\label{t:orlicz}
Let $X=L_\Phi(\mu)$ as above, and let $T\in L(X^*)$. Then
$T\in(Z^*(X))^{\com}$ if and only if
\begin{equation}\label{e:triangular}
T=AP+BQ+RP,
\end{equation}
where $A,B\in Z(X^*)$ and $R:L_\Psi\to S$ is a module map. Every module
map $L_\Psi\to S$ vanishes on $(L_\Psi)_a$, and $T\in Z(X^*)$ if and only if
$R=0$. Consequently $Z(X^*)=(Z^*(X))^{\com}$ if and only if every module
map $L_\Psi\to S$ is zero; in particular, equality holds whenever
$(L_\Psi)_a=L_\Psi$.
\end{theorem}

\begin{proof}
\emph{Sufficiency.} $AP$ and $BQ$ are central, hence commute with
$Z^*(X)\subseteq Z(X^*)$. For $g\in L^\infty$, $P$ commutes with $M_g^*$, so
$RP(g\varphi)=R(gP\varphi)=gRP\varphi$.

\emph{Necessity.} Let $T\in(Z^*(X))^{\com}$. Since $P,Q\in Z(X^*)$ commute
with every $M_g^*$, the four maps $PT|_{L_\Psi}$, $QT|_{L_\Psi}$,
$PT|_S$, $QT|_S$ are module maps.

(i) $PT|_S=0$. Let $0\le\sigma\in S$, $k=PT\sigma\in L_\Psi$, let $B_m$
be as in Lemma~\ref{l:Ssupport}, and put $C_m=\bigcap_{j\ge m}(\Omega
\setminus B_j)$. Since $1_{B_j}\to0$ a.e., $C_m\uparrow\Omega$ up to a null
set. As $0\le1_{C_m}\sigma\le1_{\Omega\setminus B_m}\sigma$, we have
$\|1_{C_m}\sigma\|\le2^{-m}$, and $1_{C_m}k=PT(1_{C_m}\sigma)$ gives
$\|1_{C_m}k\|\le\|T\|2^{-m}$. Since $1_{C_m}|k|\uparrow|k|$ a.e.\ and the
norm of $L_\Psi$ has the Fatou property, $\|k\|=\lim_m\|1_{C_m}k\|=0$.

(ii) Every module map $\Theta:L_\Psi\to S$ vanishes on $(L_\Psi)_a$; in
particular so does $QT|_{L_\Psi}$. Let $h\in(L_\Psi)_a$, $\tau=\Theta h\in
S$, and let $B_m$ be given by Lemma~\ref{l:Ssupport} for $|\tau|$. Then
$1_{B_m}\tau=\Theta(1_{B_m}h)$, and
\[
\|\tau\|\le\|1_{B_m}\tau\|+\|1_{\Omega\setminus B_m}|\tau|\|\le
\|\Theta\|\,\|1_{B_m}h\|+2^{-m}.
\]
Since $|1_{B_m}h|\le|h|$, $1_{B_m}h\to0$ a.e.\ and $h$ has order continuous
norm, $\|1_{B_m}h\|\to0$. Hence $\tau=0$.

(iii) $PT|_{L_\Psi}$ is central on $L_\Psi$: it commutes with every band
projection $k\mapsto1_Ek$ of $L_\Psi$, so Lemma~\ref{l:bandcomm} applies in
$E=L_\Psi$.

(iv) $QT|_S$ is central on $S$. Otherwise, as in the proof of
Lemma~\ref{l:bandcomm}, it is not band preserving: there is $\varphi\in S$
with $QT\varphi\notin\{\varphi\}^{dd}$. Let $P_\varphi$ be the band
projection onto $\{\varphi\}^{dd}$. It is central, so it commutes with the
$M_g^*$, and $(I-P_\varphi)QT|_{\{\varphi\}^{dd}}$ is a nonzero module map
from $\{\varphi\}^{dd}$ to the disjoint band $S\cap\{\varphi\}^d$. This
contradicts Lemma~\ref{l:Smodule}.

Put $A=PTP$, $B=QTQ$ and $R=QT|_{L_\Psi}$. By (iii) and (iv), $A$ and $B$
are central on $X^*$ (a central operator on each of two complementary bands
gives a central operator on $X^*$). By (i), $PTQ=0$, so
$T=PTP+QTP+QTQ=AP+RP+BQ$, and $R$ vanishes on $(L_\Psi)_a$ by (ii).

\emph{Centrality.} If $R=0$, then $T=AP+BQ$ is central. Conversely, if $T$
is central, then so is $RP=T-AP-BQ$. For $h\in L_\Psi$ we have $RPh=Rh\in S$,
which is disjoint from $h$, so $Rh=0$ by Lemma~\ref{l:bandpres}(a).
\end{proof}

\subsection{Equality for all Orlicz spaces}

Let $\varphi$ be the right derivative of $\Phi$ and $\psi(v)=\sup\{s\ge0:\
\varphi(s)\le v\}$ its right inverse, so that $\Phi(u)=\int_0^u\varphi$ and
$\Psi(v)=\int_0^v\psi$~\cite[Chapter~I, \S2]{KR}. Since $\varphi$ and $\psi$
are nondecreasing,
\begin{equation}\label{e:youngpt}
\Phi(u)\le u\varphi(u),\qquad\Psi(v)\le v\psi(v)\qquad(u,v\ge0),
\end{equation}
and if $v<\varphi(u)$, then $\psi(v)\le u$.

\begin{lemma}\label{l:split}
Let $x\in L_\Phi$ and $h\in L_\Psi$. There is $A\in\Sigma$ such that
$1_Ax\in E_\Phi$ and $1_{\Omega\setminus A}h\in E_\Psi$.
\end{lemma}

\begin{proof}
We may assume $x,h\ge0$. By H\"older's inequality, $xh\in L^1(\mu)$.

\emph{A weight.} Let $F_1\subseteq F_2\subseteq\cdots$ be sets of finite
measure with union $\Omega$, and $B_n=F_n\cap\{x\le n,\ h\le n\}$. Then
$B_n\uparrow\Omega$ up to a null set, and $1_{B_n}x\in E_\Phi$,
$1_{B_n}h\in E_\Psi$, since bounded functions supported on sets of finite
measure belong to $E_\Phi$ and $E_\Psi$. By dominated convergence there are
$n_1<n_2<\cdots$ with $\int_{\Omega\setminus B_{n_j}}xh\,d\mu\le4^{-j}$.
Put $w=1$ on $B_{n_1}$, $w=j$ on $B_{n_j}\setminus B_{n_{j-1}}$ for $j\ge2$.
Then $w\ge1$ is measurable,
\[
\int w^2xh\,d\mu\le\int xh\,d\mu+\sum_{j\ge2}j^2\,4^{-(j-1)}<\infty,
\]
and $\{w<\lambda\}\subseteq B_{n_j}$ whenever $j\ge\lambda$.

\emph{The set.} Let $A=\{wh\ge\varphi(wx)\}$. On $A$, by \eqref{e:youngpt},
$\Phi(wx)\le wx\,\varphi(wx)\le w^2xh$. On $\Omega\setminus A$ we have
$wh<\varphi(wx)$, hence $\psi(wh)\le wx$ and, by \eqref{e:youngpt},
$\Psi(wh)\le wh\,\psi(wh)\le w^2xh$.

\emph{Conclusion.} Let $\lambda\ge1$ and fix an integer $j\ge\lambda$. On
$A\cap\{w\ge\lambda\}$ we have $\Phi(\lambda x)\le\Phi(wx)\le w^2xh$, and
$\{w<\lambda\}\subseteq B_{n_j}$, so
\[
\int_A\Phi(\lambda x)\,d\mu\le\int w^2xh\,d\mu+
\int\Phi(\lambda1_{B_{n_j}}x)\,d\mu<\infty .
\]
Thus $1_Ax\in E_\Phi$. In the same way, using the estimate for $\Psi(wh)$ on
$\Omega\setminus A$, $\int_{\Omega\setminus A}\Psi(\lambda h)\,d\mu<\infty$
for every $\lambda>0$, i.e.\ $1_{\Omega\setminus A}h\in E_\Psi$.
\end{proof}

\begin{remark}
The proof uses only that $xh\in L^1(\mu)$. Pointwise it says that
$\min\bigl(\Phi(wx),\Psi(wh)\bigr)\le w^2xh$, a form of Young's
inequality; the point is that the scale $w$ may vary from point to point
and tend to $\infty$ while $w^2xh$ remains integrable.
\end{remark}

\begin{theorem}\label{t:orliczall}
Let $\mu$ be a $\sigma$-finite measure, $\Phi$ an $N$-function and
$X=L_\Phi(\mu)$. Every bounded module map $L_\Psi\to S$ is zero.
Consequently $Z(X^*)=(Z^*(X))^{\com}$, and every operator in
$(Z^*(X))^{\com}$ has the form $AP+BQ$ with $A,B\in Z(X^*)$.
\end{theorem}

\begin{proof}
Let $R:L_\Psi\to S$ be a bounded module map, $h\in L_\Psi$ and $\tau=Rh$. Let
$x\in X$, and let $A$ be given by Lemma~\ref{l:split}. Since
$1_{\Omega\setminus A}h\in E_\Psi=(L_\Psi)_a$, part~(ii) of the proof of
Theorem~\ref{t:orlicz} gives $1_{\Omega\setminus A}\tau=
R(1_{\Omega\setminus A}h)=0$, so $\tau=1_A\tau$. Since $1_Ax\in E_\Phi$ and
$\tau\in S=(E_\Phi)^\perp$,
\[
\tau(x)=(1_A\tau)(x)=\tau(1_Ax)=0 .
\]
As $x$ was arbitrary, $\tau=0$, so $R=0$. The remaining assertions follow
from Theorem~\ref{t:orlicz}.
\end{proof}

For the dual of an Orlicz space we need the converse of
Lemma~\ref{l:sum}.

\begin{lemma}\label{l:sumconv}
Let $X=X_1\oplus X_2$ be as in Lemma~\ref{l:sum}. If
$Z(X^*)=(Z^*(X))^{\com}$, then $Z(X_1^*)=(Z^*(X_1))^{\com}$.
\end{lemma}

\begin{proof}
Let $T_1\in(Z^*(X_1))^{\com}$ and $S\in Z(X)$. Since $S$ commutes with the
central coordinate projections, $S=S_1\oplus S_2$, and $S_i\in Z(X_i)$
because $|S_ix_i|\le\kappa|x_i|$. Hence $S^*=S_1^*\oplus S_2^*$ commutes
with $T_1\oplus0$. So $T_1\oplus0\in(Z^*(X))^{\com}=Z(X^*)$, and
$T_1\in Z(X_1^*)$.
\end{proof}

\begin{corollary}\label{c:dual}
Let $\mu$ be $\sigma$-finite and $\Phi$ an $N$-function. Then
$Z(X^*)=(Z^*(X))^{\com}$ for $X=L_\Phi(\mu)^*$.
\end{corollary}

\begin{proof}
By (O1), $X=L_\Psi\oplus S$, a direct sum of two bands with a lattice norm.
The band $S$ is an AL-space (Lemma~\ref{l:SAL}), so it has order
continuous norm and satisfies equality by Proposition~\ref{p:oc}. The space
$L_\Psi$ carries the Orlicz norm, which is equivalent to its Luxemburg norm;
passing to an equivalent lattice norm changes neither the centers nor the
commutant. With the Luxemburg norm, $L_\Psi$ satisfies equality by
Theorem~\ref{t:orliczall} applied to $\Psi$, whose complementary function is
$\Phi$. By Lemma~\ref{l:sum}, equality holds for $X$.
\end{proof}

Together, Lemmas~\ref{l:sum} and~\ref{l:sumconv} show that equality holds
for a direct sum if and only if it holds for each summand. In particular,
equality holds for $L_\Phi^*$ if and only if it holds for $L_\Psi$.

\subsection{Orlicz spaces not covered by earlier results}

Theorem~\ref{t:orliczall} covers spaces to which
Propositions~\ref{p:ideal} and~\ref{p:oc} and part~(1)
of~\cite[Theorem~9.11]{AAK} do not apply. For example, let $X=L_{\exp}[0,1]$,
$\Phi(u)=e^u-1-u$. Proposition~\ref{p:oc} does not apply, since
$\Phi\notin\Delta_2$. Nor does part~(1)
of~\cite[Theorem~9.11]{AAK}: if $\mu$ is finite and nonatomic and
$L_\Phi\ne E_\Phi$, then $L_\Phi$ has no topological order unit. Indeed,
let $0\le w\in L_\Phi$ and choose $k$ with $G=\{w\le k\}$ of positive
measure. Since $\mu$ is finite and nonatomic, $L_\Phi\ne E_\Phi$ means that $\Phi$
fails $\Delta_2$ for large $u$, and then $L_\Phi(G)\ne E_\Phi(G)$ for every
$G$ of positive measure; so there is $0\le x\in L_\Phi$
supported in $G$ with $\theta_\Phi(x)=1$. If $|y|\le nw$, then $|y|\le nk$
on $G$, so $|x-y|\ge(x-nk)^+$ on $G$ and
\[
\|x-y\|\ge\theta_\Phi\bigl((x-nk)^+\bigr)=\theta_\Phi(x)=1,
\]
because $x-(x-nk)^+$ is bounded and supported in $G$, hence in $E_\Phi$. So
$x$ is not in the closed ideal generated by $w$.

Proposition~\ref{p:ideal} does not apply either. This follows from the
preceding paragraph and the next lemma, since the Luxemburg norm has the
Fatou property in the sense of the lemma (cf.~\cite[Theorem~131.7]{Za}): if
$0\le x_k\uparrow x$ a.e.\ and $\|x_k\|\le M$, then $\int\Phi(x_k/M)\,d\mu\le1$,
so $\int\Phi(x/M)\,d\mu\le1$ by monotone convergence, i.e.\ $x\in L_\Phi$ and
$\|x\|\le M$.

\begin{lemma}\label{l:fatouqi}
Let $\mu$ be a finite measure and let $X$ be an ideal of $L^0(\mu)$ with
$L^\infty(\mu)\subseteq X$, equipped with a complete lattice norm with the
Fatou property: if $0\le x_k\uparrow x$ a.e., $x_k\in X$ and
$\sup_k\|x_k\|<\infty$, then $x\in X$ and $\|x\|\le\sup_k\|x_k\|$. If $X$
is lattice isomorphic to a closed ideal of a Banach lattice with a
topological order unit, then $X$ has a topological order unit.
\end{lemma}

\begin{proof}
Let $V$ be a lattice isomorphism of $X$ onto a closed ideal $J$ of a Banach
lattice $Y$, and let $u$ be a topological order unit of $Y$. Positive
operators between Banach lattices are bounded, so $c=\|V^{-1}\|<\infty$.
Lattice operations in $X$ are pointwise a.e. For $x\in X_+$ and $r>0$ the
element $Vx\wedge ru$ lies in $J$, since $0\le Vx\wedge ru\le Vx$.

\emph{Step 1.} Put $g_k=V^{-1}(V(k1)\wedge u)\in X_+$, $k\in\N$. For $m\ge
k$, $V(g_m\wedge k1)=V(m1)\wedge u\wedge V(k1)=V(k1)\wedge u=Vg_k$, so
$g_m\wedge k1=g_k$. Hence $(g_k)$ is increasing, and with $g=\sup_kg_k$
(pointwise, with values in $[0,\infty]$) we have $g_k=\min(g,k)$ a.e.

\emph{Step 2: $V^{-1}(Vx\wedge u)=x\wedge g$ for $x\in X_+$.} Let
$x_k=x\wedge k1$. Since $x_k\le k1$,
\[
Vx_k\wedge u=Vx_k\wedge V(k1)\wedge u=V(x_k\wedge g_k)=V(x_k\wedge g),
\]
because $x_k\wedge\min(g,k)=x_k\wedge g$. By Birkhoff's inequality in $Y$,
$0\le Vx\wedge u-Vx_k\wedge u\le Vx-Vx_k$; applying $V^{-1}$,
$0\le V^{-1}(Vx\wedge u)-x_k\wedge g\le x-x_k$, and $x-x_k\to0$ a.e. Letting
$k\to\infty$ gives the claim. Replacing $x$ by $x/n$, also
$V^{-1}(Vx\wedge nu)=x\wedge ng$ for $n\in\N$.

\emph{Step 3: $g\in X$.} By Step 1, $\|\min(g,k)\|=\|g_k\|=\|V^{-1}(V(k1)\wedge
u)\|\le c\|u\|$ for every $k$, and $\min(g,k)\uparrow g$. Moreover
$g<\infty$ a.e.: if $A=\{g=\infty\}$, then $k1_A\le\min(g,k)$, so
$k\|1_A\|\le c\|u\|$ for all $k$, whence $1_A=0$ in $X$ and $\mu(A)=0$. By the
Fatou property, $g\in X$.

\emph{Step 4.} Let $x\in X_+$ and $\varepsilon>0$. Since $u$ is a topological
order unit, there is $y\in Y$ with $|y|\le mu$ for some $m$ and
$\|Vx-y\|<\varepsilon$. For $n\ge m$,
\[
0\le Vx-Vx\wedge nu=(Vx-nu)^+\le(Vx-|y|)^+\le(Vx-y)^+\le|Vx-y|,
\]
so $\|x-x\wedge ng\|=\|V^{-1}(Vx-Vx\wedge nu)\|\le c\varepsilon$ by Step 2.
Since $0\le x\wedge ng\le ng$, every $x\in X_+$, and hence every $x\in X$, lies
in the closed ideal of $X$ generated by $g$; that is, $g$ is a topological
order unit of $X$.
\end{proof}

Thus $L_{\exp}[0,1]$, and more generally every $L_\Phi(\mu)$ with $\mu$
finite and nonatomic and $L_\Phi\ne E_\Phi$, is not lattice isomorphic to a
closed ideal of a Banach lattice with a topological order unit, and none of
Propositions~\ref{p:ideal} and~\ref{p:oc} and part~(1)
of~\cite[Theorem~9.11]{AAK} applies to it, while
Theorem~\ref{t:orliczall} applies to all of them.

\section{Open problems}\label{s:problems}

\begin{problem}
Characterize the Banach lattices $X$ (with topologically full center) for
which $Z(X^*)=(Z^*(X))^{\com}$. Is the converse of Theorem~\ref{t:D} true?
\end{problem}

\begin{problem}
For which concave functions $\psi$ does equality fail for the Marcinkiewicz
spaces $m_\psi$ and $M_\psi$? Our construction uses two disjoint positive
symmetric singular functionals living on disjoint sets of scales; the
choice $\psi(t)=\log(1+t)$ enters through $\psi(2t)/\psi(t)\to1$ in
Proposition~\ref{p:dixmier} and the scale estimates in Lemma~\ref{l:beta}
and Proposition~\ref{p:disjoint}. Note also that our proof of equality for
Orlicz spaces (Theorem~\ref{t:orliczall}) uses the fact that the singular
part of the dual of an Orlicz space is an AL-space (Lemma~\ref{l:SAL}; see
also~\cite{An,Ra,dJ}), while for some Marcinkiewicz function spaces it is
not~\cite[Examples~8.2--8.4]{dJ}.

The method needs nonzero positive symmetric functionals, and these do not
exist when $\psi$ grows regularly. Suppose that $\psi(2t)\ge(1+\delta)
\psi(t)$ for all $t>0$, with $\delta>0$. Let $D_mx(t)=x(t/m)$. Since
$(D_mx)^*(t)=x^*(t/m)$,
\[
\int_0^t(D_{2^n}x)^*=2^n\int_0^{t/2^n}x^*\le2^n\psi(t/2^n)\|x\|\le
\Bigl(\frac2{1+\delta}\Bigr)^n\psi(t)\|x\|,
\]
so $\|D_{2^n}\|\le2^{nr}$ with $r=1-\log_2(1+\delta)<1$. On the other hand,
$D_2x$ is equimeasurable with the sum of two disjoint copies of $x$, so a
positive symmetric functional $s$ satisfies $s(D_{2^n}x)=2^ns(x)$ for
$x\ge0$. Hence $2^ns(x)\le2^{nr}\|s\|\|x\|$ for all $n$, and $s=0$. The same
argument works for $m_\psi$, with $(D_mx)_n=x_{\lceil n/m\rceil}$. This
applies, e.g., to $\psi(t)=t^p$, $0<p<1$ (then $r=1-p$), i.e.\ to the weak
$L^q$ spaces, $q=1/(1-p)$. In the other direction, by~\cite[Theorem~3.4]{DPSS},
$M_\psi(0,\infty)$ carries a nonzero symmetric functional supported at
$\infty$ (i.e.\ vanishing on $M_\psi\cap L^1$) if and only if
$\liminf_{t\to\infty}\psi(2t)/\psi(t)=1$ and $\psi(t)\to\infty$ as
$t\to\infty$, and one supported at $0$ (vanishing on $M_\psi\cap L^\infty$)
if and only if $\liminf_{t\downarrow0}\psi(2t)/\psi(t)=1$ and
$\psi(0+)=0$. Such functionals are singular. In~\cite{DPSS} a positive
functional is called symmetric if it is monotone with respect to
submajorization; such a functional takes equal values on equimeasurable
positive functions, so it is symmetric in our sense. For
$\psi(t)=\log(1+t)$ only the behaviour at $\infty$ is relevant, since
$\psi(2t)/\psi(t)\to2$ as $t\downarrow0$.

The dichotomy is governed by the regularity of the function generating the
Marcinkiewicz space, studied by Mekler~\cite{Me00,Me18}. For Marcinkiewicz
spaces on $[0,1]$, with $\psi(t)=\int_0^tf$, Mekler proved that
$\liminf_{t\to0}\psi(2t)/\psi(t)>1$ is equivalent to the regularity
$f^{**}\le Kf^*$ of $f$, to the equality of the principal rearrangement
invariant and majorant ideals $N_f=M_f$, and to the normability of
$N_f$~\cite[Theorem~3]{Me00} (the main implications go back
to~\cite{Me75}; see also~\cite[Corollary~4.1]{BM77}, where these conditions
are shown to be equivalent to the Hardy--Littlewood property of $M_f$ and
to $\lim_{\tau\to\infty}\tau^{-1}\|\sigma_\tau\|_{M_f\to M_f}=0$ for the
dilations $(\sigma_\tau x)(t)=x(t/\tau)$, $\tau\ge1$), and that it is also
equivalent to a power estimate $\|\rho_s\|_{M_\psi\to M_\psi}\le ks^{-p}$,
$0<s\le1$, with $k>0$ and $p\in(0,1)$, for the compression--dilation
operators $(\rho_sx)(t)=x(st)$ if $st\le1$, $(\rho_sx)(t)=0$ if $st>1$, and to
the averaging of $N_f$ by every $\sigma$-subalgebra~\cite[Theorem~1.4]{Me18}. This suggests two
questions. \emph{Does equality hold for every Marcinkiewicz space with
regular $\psi$? Does it fail for every $m_\psi$ and $M_\psi(0,\infty)$ with
$\liminf_{t\to\infty}\psi(2t)/\psi(t)=1$ and $\psi(\infty)=\infty$, and for
every $M_\psi(0,\infty)$ with $\liminf_{t\downarrow0}\psi(2t)/\psi(t)=1$
and $\psi(0+)=0$?}
\end{problem}

\subsection*{Acknowledgment}
The counterexample and several of the other results of this paper were
obtained in an extended interaction between the second author and Claude, an
AI system developed by Anthropic. Claude also prepared drafts of the text.
Separate instances of Claude, given only the manuscript, reviewed the
proofs. The authors have checked all arguments and are fully responsible for
the contents of the paper.

\end{document}